\documentclass{article}

\usepackage{eurosym}
\usepackage{amssymb}
\usepackage{amsmath}
\usepackage{amsfonts}
\usepackage{version}
\usepackage{appendix}
\usepackage{xcolor}
\usepackage{soul}
\usepackage{cancel}
\usepackage{ulem}
\usepackage{graphicx}
\usepackage{float}
\usepackage{subcaption}
\usepackage{url}
\usepackage[english]{babel}
\usepackage{enumitem}
\usepackage{varwidth}
\usepackage{tasks}
\usepackage{booktabs,ctable,threeparttable}

\providecommand{\U}[1]{\protect \rule{.1in}{.1in}}
\newtheorem{theorem}{Theorem}

\newtheorem{condition}{Condition}

\newtheorem{corollary}{Corollary}
\newtheorem{definition}{Definition}
\newtheorem{example}{Example}

\newtheorem{notation}{Notation}
\newtheorem{proposition}{Proposition}
\newtheorem{remark}{Remark}
\newcommand{\integers}{\mathbb{Z}}
\newcommand{\autocov}{\mathbb{C}ov}
\newcommand{\var}{\mathbb{V}ar}
\newcommand{\expectation}{\mathbb{E}}

\newcommand{\reals}{\mathbb{R}}

\def\intv#1[#2..#3]{\mathopen{#1[}#2\mathrel{{.}\,{.}}\nobreak#3\mathclose{#1]}}

\begin{document}

\title{\textbf{A unified theory for ARMA models with\\
varying coefficients: One solution fits all\footnote{The proofs of the statements are included in the appendices  of this paper and they are available to the reader upon request.}}}
\author{M. Karanasos$^{\dagger \triangleright }$, A. Paraskevopoulos$%
^{\ddagger }$, T. Magdalinos$^{\ast }$, A. Canepa$^{\star }$ \\
$^{\dagger }$Brunel University London, $^{\ddagger }$University of Piraeus,\\ $^{\ast }$University of Southampton, $^{\star }$ University of Turin}
\maketitle

\begin{abstract}
For the large family of ARMA models with variable coefficients (TV-ARMA), either deterministic or stochastic, we provide an explicit and computationally tractable representation based on the general solution of the associated linear difference equation. Analogous representations are established for the fundamental properties of such processes, including the Wold-Cram\'{e}r decomposition and their covariance structure as well as explicit optimal linear forecasts based on a finite set of past observations. These results are grounded on the principal determinant, that is a banded Hessenbergian representation of a restriction of the Green function involved in the solution of the linear difference equation associated with TV-ARMA models, built up solely of the autoregressive coefficients of the model. The $L_2$ convergence properties of the model are a consequence of the absolute summability of the aforementioned Hessenbergian representation, which is in line with the asymptotic stability and efficiency of such processes. The  invertibility of the model is also a consequence of an analogous condition, but now the  Green function is built up of the moving average coefficients. The structural asymmetry between constant and deterministically time-varying coefficient models, that is the backward and forward asymptotic efficiency differ in an essential manner, is formally demonstrated. An alternative approach to the Hessenbergian solution representation is described by an equivalent procedure for manipulating time-varying polynomials. The practical significance of the theoretical results in this work is illustrated with an application to U.S. inflation data. The main finding is that inflation persistence increased after 1976, whereas from 1986 onwards the persistence declines and stabilizes to even lower levels than the pre-1976 period. \vspace{0.3in}

\vspace{0.1in}\newline
 \vspace{1in}

$^{\triangleright }${\footnotesize 
We would like to thank S. Dafnos, who contributed to the interpretation and synthesis of the work, which led to an earlier version of the paper, and P. Koutroumpis, who provided the empirical Section for this earlier version of the paper.}

\vspace{0.2cm}\newpage

{\footnotesize We gratefully} {\footnotesize acknowledge the helpful
conversations we had with L. Giraitis, and G. Kapetanios in the preparation
of the paper. We would also like to thank R. Baillie, M. Brennan, L.
Bauwens, A. Demos, W. Distaso, D. van Dijk, C. Francq, P. Fryzlewicz, E.
Guerre, C. Gourieroux, M. Guidolin, A. Harvey, C. Hommes, E. Kyriazidou, S.
Leybourne, P. Marsh, P. Minford, C. Robotti, W. Semmler, R. Smith, T. Ter%
\"{a}svirta, E. Tzavalis, P. Zaffaroni, and J-M Zakoian for suggestions and
comments which greatly improved many aspects of the paper. We are grateful
to seminar participants at Aarhus University (Creates), CREST, London School
of Economics (Department of Statistics), Erasmus University (Econometric
Institute), Imperial College (Tanaka Business School), Queen Mary University
of London (School of Economics and Finance), University of Essex (Business
School), BI Norwegian Business School, University of Copenhagen (Department
of Economics), Birkbeck College, University of London (Department of
Economics, Mathematics and Statistics), University of Nottingham (Granger
Centre for Time Series Econometrics), University of Manchester (Business
School), Cardiff University (Business School), Lancaster University
(Management School), University of Reading (ICMA Centre), King's College
London (Business School), University of Turin (Department of Economics and
Statistics), Athens University of Economics and Business (Department of
Economics) and University of Piraeus (Department of Economics). }

{\footnotesize We have also benefited from the comments given by
participants at the: SNDE 27th, 26th, 25th, 24th, 22nd and 21st Annual
Symposiums (Federal Reserve Bank of Dallas, March 2019; Keio University,
Tokyo, March 2018; ESSEC Business School, Paris, March 2017; University of
Alabama, March 2016; The City University of New York, April 2014; University
of Milano-Bicocca, March 2013, respectively), 5th and 4th AMEF Conferences
(University of Macedonia, April 2019, 2018, respectively), 2nd DSSR
Conference (University of Milano-Bicocca, February 2019), 5th, 3rd, 2nd and
1st IAAE Conferences (University of Montreal, June 2018; University of
Milano-Bicocca, June 2016; University of Macedonia, Greece, June 2015; Queen
Mary University of London, June 2014, respectively), 12th, 11th, 8th and 7th
CFE Conferences (University of Pisa, December 2018, December 2014; Senate
House, University of London, December 2013, December 2017), 17th, 13th and
12th CRETE Conferences (Tinos, Crete and Naxos, July 2018, 2014 and 2013,
respectively), SEM Conference (MIT, Boston, July 2017), 69th and 67th ESEM
European Meetings (University of Geneva, August 2016; Toulouse, France,
August 2014, respectively), Final and 1st RastaNews Annual Conferences
(Catholic University, Milan, January 2016; Palazzo Clerici, Milan, January
2014, respectively), 35th ISF Conference (Erasmus University, Rotterdam,
July 2014), 2nd ISNPS Conference (Cadiz, Spain, June 2014), NSVSM Conference
(University of Paderborn, July 2014), 9th and 8th BMRC-QASS conferences on
Macro and Financial Economics (Brunel University London, May 2014 and 2013,
respectively), and 3rd Humboldt-Copenhagen Conference on Financial
Econometrics (Humboldt University, Berlin, March 2013).}

\end{abstract}

\newpage

\section{Introduction}
Modelling time series processes with variable  coefficients has received considerable attention in recent years in the wake of several financial crises and high volatility due to frequent abrupt changes in the market. Justification for the use of such structures can be found in Timmermann and van Dijk (2013); for example, for the dynamic econometric modelling and forecasting in the presence of instability see the papers in the corresponding Journal of Econometrics special issue, i.e., Pesaran et al. (2013). For some more recent references see, for example, Cavaliere and Taylor (2005, 2008),  Cavaliere and Georgiev (2008), Giraitis et al. (2014), Harvey et  al. (2018), Chambers and Taylor (2020). Time-varying coefficient models are extensively applied by practitioners, and their importance is widely recognized (see, for example, Granger, 2007 and 2008, Petrova, 2019, Kapetanios et al., 2019, 2020).\footnote{A growing empirical literature in macroeconomics is testimony to their importance. See, for example, Evans and Honkapohja (2001, 2009).} Crucial advances in both the theory and the empirics for these structures are the works by Whittle (1965), Abdrabbo and Priestley (1967), Rao (1970), Hallin (1979, 1986), Kowalski and Szynal (1990, 1991) and Grillenzoni (1993, 2000).\footnote{See also Francq and Gautier (2004a, 2004b). We refer to the introduction of Azrak and M\'{e}lard (2006) and Alj et al. (2017) for further references.}

This paper provides a general framework for the study of autoregressive moving average models with variable coefficients and heteroscedastic errors (hereafter, TV-ARMA). There are two large classes of
stochastic processes: the ones with deterministically and those with stochastically time-varying coefficients. Both types have been widely applied in many fields of research, such as economics, finance and engineering (see for example Barnett et al., 2014, Aastaveit et al., 2017, Mumtaz and Petrova, 2018), but traditionally they have been examined separately.
The new framework unifies them by showing that one solution fits all. More specifically, we obtain explicit and computationally feasible solution  representations that generate the fundamental properties of these models, whereas the scope of the useful tool which is traditionally used to obtain such representations, that is the characteristic polynomials, is diminished when time parameter  variation is present (see, for details, Hallin, 1979, and Grillenzoni, 1990).

Following Miller (1968), in a series of papers, Hallin (1979, 1984, 1986, 1987), Singh and Peiris (1987), and Kowalski and Szynal  (1990, 1991) employ the one-sided Green's function (Green function for short) for the Wold decomposition of TV-ARMA models and their core properties. This is grounded on a particular solution of the associated linear difference equation.  However,  the Green function involved in this solution, is defined implicitly via a fundamental set of solutions, the elements of which was not, in the general case, explicitly specified. As an alternative, recursive methods were employed to compute it (see, for example, Grillenzoni, 2000, and Azrak and Melard, 2006). An explicit representation of the Green function  depends upon the availability of a fundamental set of solutions (marked also by Hallin, 1979) whose elements (known as fundamental or linearly independent solutions) must be explicitly expressed and easily handled, being an ongoing research issue. Having at our disposal such a fundamental solution set, we also have the general homogeneous and in turn the general nonhomogeneous solution of the  TV-LDE associated with TV-ARMA models.

In this work we provide such a fundamental solution set  yielding a banded Hessenbergian (determinant of a banded Hessenberg matrix) representation of the Green
function restriction involved in the solution of the associated TV-LDE (see Paraskevopoulos and Karanasos, 2021), termed here as ``principal determinant" and denoted by $\xi(t,s)$ (see  eq. (\ref{ksi})). The entries of $\xi(t,s)$ are the autoregressive coefficients taken at consecutive time instances. The first fundamental solution is the principal determinant. The remaining fundamental solutions are also expressed in terms of the principal determinant and therefore as banded Hessenbergians, yielding an explicit form for the entries of the companion matrix product. The advantages of the principal determinant representation of the Green function  are multiple. 
Using elementary properties of banded Hessenbergians, we obtain the general homogeneous solution of the linear difference equation associated with a TV-ARMA model as a linear combination of the principal determinant  times prescribed initial condition values (see eqs. (\ref{Unified homogeneous solution}, \ref{alternative homogeneous solution})). The particular solution, mentioned earlier, is also recovered, but in a more explicit form, as a linear combination of the principal determinant times consecutive instances of the forcing term (see eq. (\ref{PARTICULAR SOLUION})). These results lead to an efficient interpretation of the model, described as a decomposition into four structurally different parts (see eq. (\ref{TVAR(p)SOL1})). In this interpretation the fundamental solutions could be portrayed as autoregressive coefficients of prescribed random variables.

The Hessenbergian representation of TV-ARMA processes  gives rise to  an explicit condition, that is the absolutely summability of terms involving the principal determinant (see (\ref{abs summability condition})), which leads to easily handled explicit representations of their fundamental properties in $L_2$ spaces. The convergence to zero of the principal determinant along with the boundedness of the autoregressive coefficients ensure the asymptotic stability of all bounded solution processes in $L_2$. The absolute summability condition along with the boundedness of the drift and the moving average coefficients guarantee the existence and uniqueness of the Wold-Cram\'{e}r decomposition of asymptotically stable solution processes associated with finite unconditional first two moments and autocovariance function. The invertibility of such solution process is guaranteed by the boundedness of the autoregressive coefficients and the absolutely summability of terms generated by the principal determinant, whereas, in this case, its nonzero entries are the moving average coefficients of the model (see Section \ref{Invertibility}). As a consequence, we derive explicitly the optimal linear predictor for the asymptotically stable solution process of TV-ARMA models along with the forecasting error and its associated mean square error, when an infinite set of data is observed. An alternative approach, employs the aforementioned explicit representation of the model, which yields its optimal linear predictor when a finite set of data is observed. Therefore the problem of obtaining such a  predictor is reduced to a matter of dealing with a linear regression. 
In a causal environment, it turns out  that both approaches yield identical forecasting and mean square errors and therefore identical optimal predictors.  
Moreover, we illustrate mathematically one of the focal points in Hallin's (1986) analysis concerning the asymptotic efficiency of such models. Namely, that in a time-varying setting two forecasts with identical forecasting horizons, but at different times, yield different mean squared errors. It turns out that the backward asymptotic efficiency is, in the general case, different from the forward one 
(termed by Hallin Granger-Andersen). In this context, we explicitly provide sufficient conditions for the forward asymptotic boundedness and uniform boundedness of the mean square forecasting error. Equally importantly, we demonstrate how the linear algebra techniques, used to obtain the general solution, are equivalent to a simple procedure for manipulating polynomials with variable coefficients (see Section \ref{SECTVOPERATORS}). In order to do so we employ the expression of the principal determinant in conjunction with the so called skew multiplication operator or symbolic operator (see, for example, Hallin, 1986, and Mrad and Farag, 2002).

Various processes with stochastic coefficients are treated within our unified framework in Section \ref{SECSTV}. We investigate two general models. First, for the random coefficients model, we show that when the principal determinant converges almost surely, then the process converges in distribution. Second, the double stochastic autoregressive model is also employed to formulate explicitly some of its fundamental properties. For example, the convergence in $L_{2}$ of the principal determinant is a necessary condition for the unconditional variance to exist.

Banded Hessenbergians are computationally tractable due to the linear running time for their calculation (see Paraskevopoulos and Karanasos, 2021). Compact solution representations of banded
Hessenbergians, established in Marrero and Tomeo (2012, 2017) and Paraskevopoulos and Karanasos (2021), can be applied to derive analogous explicit representations for the principal determinant. These results modernize and
enhance the explicit representations of time-varying models and their fundamental properties, by compact representations.

The definition of the principal determinant  arises naturally from Paraskevopoulos (2014), who (by introducing the so called infinite Gauss-Jordan elimination algorithm) provides banded Hessenbergian representations for the fundamental solution set associated with a TV-ARMA($p,q$) model (see Paraskevopoulos and Karanasos 2021). The results of the present paper are established independently of the infinite elimination algorithm.

To sum up, the unified theory enable us to study linear time series models with either stochastically or deterministically varying coefficients, thus allowing us to make theoretical comparisons between these two large classes of models.

This paper concludes with an empirical application on inflation persistence in the United States (see Section \ref{sec. InflationData}), which employs a time-varying model of inflation dynamics grounded on statistical theory. In particular, we estimate an autoregressive process with abrupt structural breaks and we compute an alternative measure of second-order time dependent persistence, which distinguishes between changes in the dynamics of inflation and its volatility and their persistence. Our main conclusion is that persistence increased after 1976, whereas from 1986 onwards it declines and stabilizes to even lower levels than the pre-1976 period. Our results are in line with those in Cogley and Sargent (2002), who find that the persistence of inflation in the United States rose in the 1970s and remained high during this decade, before starting a gradual decline from the 1980s until the early 2000s.

The outline of the paper is as follows. Section \ref{SecTVARMA} introduces the notation used in the paper followed by the principal determinant. The next Section presents the explicit representation for an extensive family of ``time-varying" ARMA models, based 
on the general solution of the associated TV-LDE. In Section \ref{sec. Stability}, we obtain a necessary and sufficient condition which guarantees the asymptotic stability of these processes. Section \ref{Representations in L_2} presents explicit formulas for the fundamental properties of the model, including the Wold-Cram\'{e}r decomposition, the unconditional moments and the autocovariance function.
Section \ref{SecSecMom} deals with the invertibility and forecasting issues. In Section \ref{SECSTV}, we examine AR models in which the drift and the autoregressive coefficients are stochastically varying. In Section \ref{SECTVOPERATORS}, we introduce a simple procedure for manipulating polynomials with variable coefficients. The next Section gives an illustrative example with abrupt structural breaks and proposes a new measure of time-varying persistence. In Section \ref{sec. InflationData} we present an empirical study for inflation persistence. The final Section of the paper contains some concluding remarks and future work.

The proofs of the statements concerning L of the current paper are demonstrated in (Paraskevopoulos and Karanasos 2021), while the proofs of stochastically oriented statements are included in the Appendix of this paper, which follows the enumeration structure of the main body of the paper. Further details are provided in the Online and the Additional Appendices.

\section{Time-varying ARMA\label{SecTVARMA}}

The aim of this Section is to provide a definition of time-varying ARMA models, hereafter termed TV-ARMA, along with the main mathematical tool for their analysis. This is the banded Hessenbergian representation of the Green function  associated with such models, which is also referred to here as the \textit{principal determinant}.

\subsection{The Model}\label{SubsecNotDef}

This Subsection introduces suitable notation and defines the basic process.
Throughout the paper we adopt the following notational conventions: The set of integers (resp. positive and non-negative integers) is denoted by $\integers$
(resp. $\integers_{>0}$ and $\integers_{\geq 0}$).  
Similarly, the set of real numbers (resp. positive and non-negative real numbers) is denoted by $\reals$  (resp. $\reals_{>0}$ and $\reals_{\geq 0}$). Moreover, 
$(\Omega ,\mathcal{F},P)$ stands for a probability space and $L_{2}(\Omega, \mathcal{F},P)$ (in short $L_{2}$) stands for the Hilbert space of random variables with finite first two moments defined on $(\Omega ,\mathcal{F},P)$. 

Let $p\in \integers_{\ge 1}$ and $q\in \integers_{\ge 1}$. A TV-ARMA($p,q$) model generates a stochastic process $\{y_t\}_{t}\in\integers$ (or simply $y_t$ process), which satisfies 
\begin{equation}\label{TVAR(P)}
y_{t}=\varphi (t)+\sum\limits_{m=1}^{p}\phi _{m}(t)y_{t-m}+u_{t},\ \text{for\ all}\ t\in\integers  
\end{equation}%
with moving average term $u_{t}$ given by \vspace{-0.15in} 
\begin{equation*}
u_{t}=\varepsilon _{t}+\sum\limits_{l=1}^{q}\theta _{l}(t)\varepsilon _{t-l},
\end{equation*}%
where the coefficients $\phi_m(t), \theta_l(t)$ are either deterministic (real valued functions) or stochastic (random variables),  $\varphi (t)$ is the time-varying drift generating a deterministic or stochastic process, while $\{\varepsilon _{t}\}_{t}$ is a mean zero random process (that is $\expectation(\varepsilon _{t})=0$) such that  $\expectation(\varepsilon _{t}\ \varepsilon _{s})=0$ for $s\not=t$ (uncorrelatedness condition), $\expectation(\varepsilon_{t}|y_{s}, s<t)=0$ for all $t$ (that is $\{\varepsilon_{t}\}$ is a martingale difference sequence relative to $\{y_{s}\}$) and the time-varying variance $\expectation(\varepsilon_{t}^{2})=\sigma ^{2}(t)$ is non-zero and bounded, that is $0<\sigma^{2}(t)\leq M<\infty $ for all $t$ and some $M\in \reals_{>0}$. The above conditions guarantee that $\varepsilon_{t}\in L_2$ and $\varepsilon_{t}\perp \varepsilon_{s}$ ($\varepsilon_{t}, \varepsilon_{s}$ are orthogonal), whenever $s\not=t$. 

If the non-constant coefficients $\phi_m(t), \theta_l(t)$ and the drift $\varphi (t)$ in eq. (\ref{TVAR(P)}) are deterministic (resp. stochastic) we shall refer to it as DTV-ARMA (resp. STV-ARMA).\footnote{Notice that in our setting the time-varying coefficients can depend on the length of the series as well, as in Azrak and M\'{e}lard  (2006); see the examples in Section \ref{Autocovariance Function}.} Further specifications and adjustments for STV-ARMA models are presented in Section \ref{SECSTV}. 

The forcing term $\upsilon_{t}$ is assigned to be the time-varying drift plus the moving average term:
\begin{equation}\label{Forcing term}
\upsilon_{t}=\varphi (t)+u_{t}.
\end{equation}
Proceeding with the notation of eq. (\ref{Forcing term}), the associated TV-LDEs($p$) of eq. (\ref{TVAR(P)}) is defined by: \begin{equation}\label{Associated TV-LDE}
y_{t}=\sum_{m=1}^{p}\phi_{m}(t)y_{t-m}+ \upsilon_t,\ \text{for\ all}\ t\in\integers, 
\end{equation}
where $\phi_{m}(t)$ stand for the autoregressive coefficients of eq. (\ref{TVAR(P)}) and  $\upsilon_t$,  $y_{t-m}$ for $1\le m\le p$  are considered as realizations of the homonyms random variables. 

In this work, both assumptions of stationarity and homoscedasticity have been relaxed (see
also, among others, Singh and Peiris, 1987, Kowalski and Szynal, 1990, 1991,
and Azrak and M\'{e}lard, 2006), which is likely to be violated in practice
and we allow $\{\varepsilon_{t}\}$ to follow, for example, a stochastic
volatility or a time-varying GARCH type of process (see, for example, the
earlier versions of the current paper: Karanasos et al., 2014c, and Karanasos et al., 2017) or we allow for abrupt
structural breaks in the variance of $\varepsilon_{t}$ (see the example in
Section \ref{Example}).

The TV-ARMA($p,q$) model nests both the TV-AR$(p)$ as a special case when $q=0$ and the ARMA$(p,q)$ specification when the drift, the autoregressive
and moving average coefficients, and the variances are all constants,
adopting for this purpose the conventional identifications: $\varphi
(t)=\varphi ,\ \phi _{m}(t)=\phi _{m},\ \theta _{l}(t)=\theta _{l},\sigma^{2}(t)=\sigma ^{2}$ for all $t$.

The relation between the process under consideration and its innovations is essentially described by the Wold-Cram\'{e}r decomposition (see Section \ref{WoldCramerDecomposition}), which is the main analytical tool for studying the  asymptotic efficiency of the model. In this case, the latest time-point of the observed random variables, denoted here by $s$, moves to the remote past ($s\rightarrow -\infty$), while the forecast time-point, denoted here by $t$, is kept fixed. The forward asymptotic efficiency of the model (so-called by Hallin, 1986, Granger-Andersen) is strongly related to the forecasting problem. It directs attention to the asymptotic properties of the mean square forecasting error (MSE for short), as the time $t$ moves to the far future, while $s$, is kept fixed (see Section \ref{FORWASYMPTEFFIC}).

One of the goals of this work is to obtain an explicit inverse of the time-varying autoregressive (AR) polynomial associated with eq. (\ref{TVAR(P)})
being denoted by $\Phi_{t}(B)$, where $B$ is the backshift or lag operator (see Section \ref{SECTVOPERATORS}). In a time-varying environment, the usual procedure employs the Green function  instead of the characteristic polynomials, which are used in the time invariant case. More specifically, we employ the principal determinant, coupled with the so called  multiplication skew operator (see the analysis in Section \ref{SUBSKEW}) to obtain the inverse of $\Phi_{t}(B)$.

We should also mention that Kowalski and Szynal (1991) used the product of companion matrices to obtain the associated Green function. Paraskevopoulos and Karanasos (2021) capitalized on the connection between the product of
companion matrices and time-varying stochastic difference equations but in the opposite direction. That is, they went the other way around and by finding an explicit and compact representation of the fundamental solutions associated with TV-ARMA models, they obtained an analogous representation for the elements of the associated companion matrix product.

\subsection{The Principal Determinant}\label{Principal Determinant}

In this and next Subsection, we provide some results on linear difference equations of order $p$ with variable coefficients (for short TV-LDEs($p$)) and their banded  Hessenbergian solution representation.

We start with the main mathematical tool of this paper, the principal determinant, denoted as $\xi (t,s)$. It has been shown in (Paraskevopoulos and Karanasos 2021), that $\xi (t,s)$ is a solution of the homogeneous linear difference equation associated with eq. (\ref{Associated TV-LDE}), that is
\begin{equation}\label{HOMOGEN}
y_{t}=\sum_{m=1}^{p}\phi_{m}(t)y_{t-m}
\end{equation}
taking on the initial conditions  $y_{s+1-p}=0,...,y_{s-1}=0, y_{s}=1$.

Amongst the various implicit representations of the  Green function in terms of undetermined fundamental sets of solutions (see Section \ref{SUBSUBHOM} and (Paraskevopoulos and Karanasos 2021) for further details), the principal determinant provides an explicit and computationally tractable representation of the Green function  restriction involved in the solution  of the associated TV-LDE (see for details Paraskevopoulos and Karanasos 2021). As a consequence, the main properties of time-varying models acquire analogous representations, as stated in Sections \ref{sec. Stability}, \ref{Representations in L_2} and \ref{SecSecMom}.

To distinguish scalars from vectors we adopt lower and uppercase boldface
symbols within square brackets for column vectors and matrices respectively: 
$\mathbf{x}=[x_{i}]$, $\mathbf{X}=[x_{ij}]$. Row vectors are indicated
within round brackets and usually appear as transpositions of column
vectors: $\mathbf{x}^{\prime}=(x_{i})$.

For every pair $(t,s)\in \integers^{2}$ such that $k=t-s\geq 1$ the principal matrix associated with the AR part of eq. (\ref{TVAR(P)}), is defined by 
\begin{equation}\label{TVAR(P)1}
\mathbf{\Phi }_{t,s}=\left[ 
\begin{array}{cccccccc}
{\footnotesize \phi }_{1}{\footnotesize (s+1)} & {\footnotesize -1} &  &  & 
&  &  &  \\ 
{\footnotesize \phi }_{2}{\footnotesize (s+2)} & {\footnotesize \phi }_{1}%
{\footnotesize (s+2)} & \ddots &  &  &  &  &  \\ 
{\footnotesize \vdots } & {\footnotesize \vdots } & \ddots & \ddots &  &  & 
&  \\ 
{\footnotesize \phi }_{p}{\footnotesize (s+p)} & {\footnotesize \phi }_{p-1}%
{\footnotesize (s+p)} & \ddots & \ddots & \ddots &  &  &  \\ 
& {\footnotesize \phi }_{p}{\footnotesize (s+p+1)} & \ddots & \ddots & \ddots
& \ddots &  &  \\ 
&  & {\small \ddots } & {\small \ddots } & {\small \ddots } & {\small \ddots 
} & \ddots &  \\ 
&  &  & {\footnotesize \phi }_{p}{\footnotesize (t-1)} & {\footnotesize \phi 
}_{p-1}{\footnotesize (t-1)} & {\small \cdots } & {\footnotesize \phi }_{1}%
{\footnotesize (t-1)} & {\footnotesize -1} \\ 
&  &  &  & {\footnotesize \phi }_{p}{\footnotesize (t)} & {\small \cdots } & 
{\footnotesize \phi }_{2}{\footnotesize (t)} & {\footnotesize \phi }_{1}%
{\footnotesize (t)}%
\end{array}%
\right],  
\end{equation}%
here and in what follows empty spaces in a matrix have to be replaced by zeros. $\mathbf{\Phi }_{t,s}$ is a lower Hessenberg matrix of order $k$. It is also a banded matrix with total bandwidth $p+1$ (the number of its non-zero diagonals, i.e., the diagonals whose elements are not all identically zero), upper bandwidth $1$ (the number of its non-zero
super-diagonals), and lower bandwidth $p-1$ (the number of its non-zero sub-diagonals). In particular, the elements of $\mathbf{\Phi }_{t,s}$ are: $(-1)$ occupying the entries of the superdiagonal, the values of the first autoregressive coefficient $\phi _{1}(\cdot )$ (from time $s+1$ to time $t$), occupying the entries of the main diagonal, the values of the $(1+r)$-th autoregressive coefficient $\phi _{1+r}(\cdot )$ for $r=1,2,\ldots ,p-1$ (from time $s+1+r$ to time $t$), occupying the entries of the $r$-th
sub-diagonal, and zero entries elsewhere. It is clear that for $p\geq k$, $\mathbf{\Phi }_{t,s}$ is a full lower Hessenberg matrix.

For every pair $(t,s)\in \integers^{2}$ with $s<t$, the principal determinant associated with eq. (\ref{TVAR(P)1}) is given by: 
\begin{equation}  \label{ksi}
\xi (t,s)=\det (\mathbf{\Phi }_{t,s}).
\end{equation}
Formally $\xi (t,s)$ is a lower Hessenbergian (determinant of a lower
Hessenberg matrix; for details on Hessenbergians see, for example, the book
by Vein and Dale, 1999). We further extend the definition of $\xi (t,s)$ so as to be defined over  $\integers^{2}$ by assigning the initial conditions: 
\begin{equation}  \label{Initial values}
\xi (t,s)=\left \{ 
\begin{array}{ccc}
1 & \text{if} & t=s \\ 
0 & \text{if} &  t<s.
\end{array}%
\right.
\end{equation}%
Under these initial values, for each
fixed $s$ the sequence $\{ \xi (t,s)\}_{t\geq s+1-p}$ turns out to be a solution sequence of eq. (\ref{HOMOGEN}) called principal fundamental sequence (see Paraskevopoulos and Karanasos 2021).
It is well known that the restriction of the Green function, often designated by $H(t,s)$ for $t>s$ and $s$ fixed, solves eq. (\ref{HOMOGEN}), assuming the prescribed values $H(s,s)=1$ for $p\ge 1$ and $H(t,s)=0$ for $t=s+1-p,...,s-1$ and $p\ge 1$ (see for example Agarwal, 2000, p.77, Property (iii) or Lakshmikantham and Trigiante, 2002, Theorem 3.4.1 p.87 ). Since the principal determinant $\xi (t,s)$ solves eq. (\ref{HOMOGEN}) under the same initial values, just as like with $H(t,s)$, the uniqueness of the solution for an initial value problem entails that for every arbitrary but fixed $s\in \integers$ the principal determinant coincides with the Green  function, that is 
$H(t,s)=\xi (t,s) $ for all $t$ such that $t\geq s+1-p$ and $p\ge 1$ (a proof from first principles is presented in Paraskevopoulos and Karanasos, 2021, Proposition 2). Therefore we will make
use of two different, although equivalent terminologies: principal determinant or Green's function with fixed $s$. An explicit and compact
representation of $\xi (t,s)$, called Leibnizian representation, is provided in the previously cited reference. An analogous representation of Hessenbergians, as a nested sum, and their inverses are established by Marrero and Tomeo (2017, 2012), respectively. 

We conclude this Section with an illustrative example, concerning the ARMA($p,q$) model with constant coefficients.

\begin{example}
\label{ex. Widom's formula} The AR polynomial $\Phi (B)=1-\sum_{m=1}^{p}\phi
_{m}B^{m}$ associated with eq. {\rm (\ref{HOMOGEN})}, whenever $\phi _{m}(t)=\phi_{m}$ (constant autoregressive coefficients), is explicitly expressed in terms of the characteristic values as $\Phi (B)=\prod_{m=1}^{p}(1-\lambda_{m}B)$. The definition in eq. {\rm (\ref{ksi})} gives the determinant of a banded Hessenbergian Toeplitz matrix (for details on Toeplitz matrices see, for example, the book by Gray,
2006)\footnote{For the use of Toeplitz matrices on double-differenced AR(1) models see Han (2007).} and satisfies the identity (usually called Widom's determinant formula; see Widom, 1958)
\begin{equation}\label{Toeplitz}
\xi _{k}=\left\vert 
\begin{array}{ccccccccc}
\phi _{1} & -1 &  &  &  &  &  &  &  \\ 
\phi _{2} & \phi _{1} & -1 &  &  &  &  &  &  \\ 
\phi _{3} & \phi _{2} & \phi _{1} & \ddots &  &  &  &  &  \\ 
\vdots & \vdots & \vdots & \ddots & \ddots &  &  &  &  \\ 
\phi _{p} & \phi _{p-1} & \phi _{p-2} & \ddots & \ddots & \ddots &  &  &  \\ 
& \phi _{p} & \phi _{p-1} & \ddots & \ddots & \ddots & \ddots &  &  \\ 
&  & \phi _{p} & \phi _{p-1} & \phi _{p-2} & \cdots & \phi _{1} & -1 &  \\ 
&  &  & \phi _{p} & \phi _{p-1} & \cdots & \phi _{2} & \phi _{1} & -1 \\ 
&  &  &  & \phi _{p} & \cdots & \phi _{3} & \phi _{2} & \phi _{1}%
\end{array}%
\right\vert ={\displaystyle\sum\limits_{m=1}^{p}}\frac{\lambda _{m}^{k+p-1}}{%
{\displaystyle\prod\nolimits_{\underset{n\neq m}{n=1}}^{p}}(\lambda
_{m}-\lambda _{n})},  
\end{equation}
where $\xi (t,s)$ is now denoted as $\xi _{k}$, since it depends only on the
forecasting horizon $k$ ($k=t-s$ indicates the order of the matrix too). The
second equality in eq. {\rm(\ref{Toeplitz})} follows (only if $\lambda_{m}\neq
\lambda _{n}$) from standard results in ARMA models (see, for example, Karanasos, 1998, 2001, and in particular, eq. (2.6) in his Corollary $2$; see also
Hamilton, 1994, pp. 12-13).
\end{example}

\section{TV-ARMA Representations}\label{SubsecMainTheor}

In the following Subsections we shall use the principal determinant to describe explicitly the general homogeneous and  particular solutions of TV-LDE($p$) in  eq. (\ref{Associated TV-LDE}). The sum of these solutions yields the general solution of eq. (\ref{Associated TV-LDE}),\footnote{In linear algebra there have been some isolated attempts to deal with the problem, which have been criticized on a number of grounds. For example, Mallik (1998) provides an explicit solution for the aforementioned equations, but it appears not to be computationally tractable (see also Mallik, 1997 and 2000). Lim and Dai (2011) point out that \textquotedblleft although explicit solutions for general linear difference equations are given by Mallik (1998), they appear to be unmotivated and no methods of solution are discussed".} which lead to an explicit and computationally tractable representation of TV-ARMA processes (see  Theorem \ref{TheoGenSol}). A useful decomposition of the innovation part of the solution is presented in Proposition \ref{PROPDECOMP}.

\subsection{Homogeneous Solution}
\label{SUBSUBHOM}
Let $a\in\integers$ and $\integers_{a}\overset{\rm def}{=}\{z\in\integers: z\ge a\}$. Let also $s,t\in \integers$ and $s<t$. Taking into account that  $s,t\in\integers_{s+1-p}$, the general solution of eq. (\ref{HOMOGEN}) over $\integers_{s+1-p}$ is denoted by $y_{t,s}^{hom}$. The latter can be expressed in terms of the principal determinant, the autoregressive coefficients and the initial conditions $y_{s+1-m}$ for $m=1,2,...,p$, as follows 
\begin{equation}\label{Unified homogeneous solution}
y_{t,s}^{hom}=\displaystyle\sum_{m=1}^{p}\sum_{r=1}^{p+1-m}\phi
_{m-1+r}(s+r)\xi (t,s+r)y_{s+1-m},  
\end{equation}%
(a proof of eq. (\ref{Unified homogeneous solution}) is provided in (Paraskevopoulos and Karanasos 2021). The two variable solution notation $y^{hom}_{t,s}$ is consistent with the analogous notation used for the principal determinant (or the Green function) and is essential for the TV-ARMA($p,q$) representation in eq. (\ref{TVAR(p)SOL1}). Moreover, the principal determinant exactly matches the specific restriction of the Green function involved in the general homogeneous solution in eq. (\ref{Unified homogeneous solution}) and, in turn, of the nonhomogeneous solution of the associated TV-LDE, that is the 
$H(t,j)=\xi (t,j) $ for all $t,j$ ($j=s+r$) such that $t\geq s+1-p$ and $s\le j\le t+p-1$ (a proof is presented in Paraskevopoulos and Karanasos, 2021, Theorem 5 and Corollary 2).

To the extent of our knowledge, the expression of the homogeneous solution in terms of the  obtained in eq. (\ref{Unified homogeneous solution}) is first recorded in earlier versions of this paper (see also Paraskevopoulos and Karanasos, 2021). 
The effects of  eq. (\ref{Unified homogeneous solution}) on this work are multiple. It makes it possible to construct the general solution of a TV-LDE($p$) entirely in terms of the Green function and any set of prescribed random variables $\{y_{s+1-p},...,y_s\}$ at a time-point $s$ of the recent or remote past. Explicit conditions guaranteeing the zero convergence of eq. (\ref{Unified homogeneous solution}) yield the asymptotic  stability and the Wold-Cram\'{e}r decomposition of the model. Another consequence, is an alternative approach to the forecasting, presented in Section \ref{SUBPREDICTOR1}.  As in the case of processes of order (1,1), this is a consequence of the explicit representation of TV-ARMA($p,q$) models in eq. (\ref{TVAR(p)SOL1}).

In all that follows the following notation is used:
\begin{equation}  \label{Exdef:xi^m1}
\xi^{(m)}(t,s)\overset{\mathrm{def}}{=}\sum_{r=1}^{p+1-m}\phi_{m-1+r}(s+r)\xi (t,s+r)
\end{equation}
(superscripts within parentheses or brackets [e.g., $(.)^{(m)}$] designate the position of the corresponding term [e.g., $m$-th term] in a
sequence, so as to distinguish position indices from power exponents). The principal determinant is identified with $\xi ^{(1)}(t,s)$, that is $\xi(t,s)\overset{\mathrm{def}}{=}\xi ^{(1)}(t,s)$. Applying eq. (\ref{Exdef:xi^m1})
with $s=t-1$, on account of $\xi (t,t)=1$ and $\xi (t,t+j)=0$ for $j\ge 1$ (see eq. (\ref{Initial values})), we conclude that: $\xi^{(m)}(t,t-1)=\phi_{m}(t)$.

It turns out (see Paraskevopoulos and Karanasos 2021) that for every $m$ with $1\le m\le p$ the sequence $\{\xi ^{(m)}(t,s)\}_{t\geq s+1-p}$ is the solution of eq. (\ref{HOMOGEN}) under the prescribed initial values: 
$y_{s+1-m}=1$ and $y_{s+1-r}=0$ for $1\leq r\leq p$ with  $r\not=m$.
Applying the expression in eq. (\ref{Exdef:xi^m1}) to the right-hand side of eq. (\ref{Unified homogeneous solution}), the homogeneous solution takes a
more condensed form: 
\begin{equation}\label{alternative homogeneous solution}
y_{t,s}^{hom}=\sum_{m=1}^{p}\xi ^{(m)}(t,s)y_{s+1-m}.
\end{equation}%
Under the initial values $y_{s}=1$ for $p\ge 1$ and $y_{s-r}=0$ for $1\le r\le p-1$, whenever $p\ge 2$, the right-hand side of eq. (\ref{alternative homogeneous solution}) turns into $\xi ^{(1)}(t,s)=\xi (t,s)$ and thus the homogeneous solution in eq. (\ref{alternative homogeneous solution})
coincides with the principal fundamental sequence $\{\xi (t,s)\}_{t\geq s+1-p}$.
In  (Paraskevopoulos and Karanasos 2021), it is shown that the set 
\begin{equation*}
\varXi_{s}=\{ \xi ^{(1)}(t,s),\xi ^{(2)}(t,s),...,\xi ^{(p)}(t,s):\ t\geq
s+1-p\}
\end{equation*}%
is a fundamental (or linearly independent) set of solutions associated with
eq. (\ref{HOMOGEN}). Moreover, it is shown there that the $m$-th fundamental solution $\xi^{(m)}(t,s) $ can be expressed as a single banded Hessenbergian too. The only difference between any two of these fundamental solutions lies in the first column (see Paraskevopoulos and Karanasos 2021).

\subsection{A Particular Solution and its  Decomposition}\label{particular solution}
A particular solution to eq. (\ref{Associated TV-LDE}) subject to the initial condition values $y_{s}=y_{s-1}=...=y_{s+1-p}=0$ is given by 
\begin{equation}  \label{PARTICULAR SOLUION}
y_{t,s}^{par}=\sum_{r=s+1}^{t}\xi (t,r)[\varphi (r)+u_{r}].
\end{equation}
A proof of the above formula is demonstrated in (Paraskevopoulos and Karanasos 2021). The solution in eq. (\ref{PARTICULAR SOLUION}) depends both on $t$ and $s$. This has to be compared with the equivalent result presented in Miller (1968, p. 40, eqs. (2.8) and (2.9)).

In Proposition \ref{PROPDECOMP} below, we introduce a decomposition of the innovation part in the particular solution in  eq. (\ref{PARTICULAR SOLUION}), which is used throughout this paper. But first we will introduce the following definition:

\begin{definition}
\label{DEFKSIQ} First, we define the function on $\integers^2$:
\begin{equation}\label{KSIq}
\xi_{q}(t,r)\overset{\mathrm{def}}{=} 
\xi(t,r)+\displaystyle\sum_{l=1}^{q}\xi(t,r+l)\theta_{l}(r+l).
\end{equation}
Second, for each $s\in\integers$, we define the function on $\integers^2$:
\begin{equation}\label{KSIsq}
\xi_{s,q}(t,r) \overset{\mathrm{def}}{=}  \displaystyle\sum_{l=s+1-r}^{q}\xi(t,r+l)\theta_{l}(r+l).
\vspace{-0.05in} 
\end{equation}
\end{definition}
As $\xi_{q}(t,r)$ is equal to $\xi (t,r)$ plus a sum of terms consisting of the first $q$ instances of $\xi(t,r+l)$  multiplied by corresponding moving average coefficients, it can also be expressed as a
banded Hessenbergian (the proof is deferred to the Online Appendix Subsection F1). The same applies to $\xi_{s,q}(t,r)$. Therefore, we shall refer to $\xi_{q}(t,r)$ and $\xi _{s,q}(t,r)$ as banded Hessenbergian coefficients. In view of eq. (\ref{KSIq}) $\xi_q(t,t)=1$, since $\xi(t,t+j)=0$, whenever $j\ge 1$ and the second term (i.e., the sum), when applied for $r=t$, vanishes. For the same reason  $\xi_{q}(t,r)=0$, whenever $r\ge t+1$. To summarize, $\xi_{q}(t,r)$ coincides with the corresponding initial conditions of $\xi(t,r)$, whenever $r\ge t$ (see eq. (\ref{Initial values})). In view of eq. (\ref{KSIsq}), if $r<s+1-q$ (or $s+1-r>q$), then the lower summation limit exceeds the upper limit, whence $\xi_{s,q}(t,r)=0$.
Finally, for a pure AR($p$) model, that is when $q=0$, we have: $\xi_{0}(t,r)=\xi (t,r)$ and $\xi_{s,0}(t,r)=0$. The banded Hessenbergian coefficient $\xi_{q}(t,r)$ must be compared with the function  defined by Peiris (1986), see his notation ``$\rm g(t,s)$" in eq. (2.2).   

\begin{proposition}
\label{PROPDECOMP} The innovation part of the
particular solution in eq. {\rm (\ref{PARTICULAR SOLUION})} can be decomposed into two parts as follows:
\begin{equation}  \label{eq. condensed}
\sum_{r=s+1}^{t}\xi (t,r)u_{r}=\sum \limits_{r=s+1}^{t}\xi
_{q}(t,r)\varepsilon _{r}+\sum \limits_{r=s+1-q}^{s}\xi
_{s,q}(t,r)\varepsilon_{r}.\vspace{-0.05in} 
\end{equation}
\end{proposition}
In the first summation of the right-hand side of eq. (\ref{eq. condensed}), the length of the time interval between $s$ and $t$ is the forecasting horizon ($k=t-s$). As a consequence, the errors whose index ranges over the integer interval $\intv\big[s+1 .. t]$~\footnote{If $a,b\in\integers$, we adopt the notation: $\intv\big[a .. b]\overset{\rm def}{=}\{z\in\integers: a\le z\le b\}$.} are unobservable. In the second one the time interval extends from $s+1-q$ to $s$ and therefore the errors whose index ranges over $\intv\big[s+1-q .. s]$ are observable.  Notice that $\xi_q(t,r)$ for $s+1\le r\le t$ can be equivalently expressed as $\xi_q(t,t-r)$ for  $0\le r\le k-1$. Similarly  $\xi_{s,q}(t,r)$ for $s+1-q\le r\le s$ can be equivalently expressed as $\xi_{s,q}(t,t-r)$ for
$k\le r\le k-1+q$ (see eq. (\ref{KSIQTS})).

Eq. (\ref{eq. condensed}) coupled with   eq. (\ref{TVAR(p)SOL1}) (in the following Subsection) will be used to obtain TV-ARMA($p,q$) forecasts in Subsections  \ref{SUBPREDICTOR1} and \ref{Double Stochastic AR Models}. If $s=t-1$ (or $k=1$), the first sum in the right-hand side of eq. (\ref{eq. condensed}) reduces to $\varepsilon_{t}$ and the second sum reduces to $\sum\nolimits_{l=1}^{q}\theta_{l}(t)\varepsilon_{t-l}$ (a formal proof is provided in the Online Appendix Subsection F2). Notice that the above mentioned reduction is in line with Remark \ref{CORk=1} below. 

\subsection{Explicit Model Representation}
\label{subsection General solution} 
The results of the previous Section on linear difference equations, are applied herein to obtain in the following Theorem an equivalent explicit representation of the stochastic process $y_t$ in eq. {\rm (\ref{TVAR(P)})} in terms of prescribed random variables: $y_s,...,y_{s+1-p}$ for any $s\in \integers$ with $s<t$. This is a consequence of the general solution representation of the  TV-LDE($p$) in  eq. (\ref{Associated TV-LDE}) as a sum of the homogeneous and the particular solutions in eqs. (\ref{alternative homogeneous solution}, \ref{PARTICULAR SOLUION}), respectively, that is $y_t=y^{hom}_{t,s}+y^{par}_{t,s}$ for all $t,s$ with $s<t$.
Notice that $y_t$ is independent of $s$, while, in context with time-series, both homogeneous and particular solution functions depend on $t$ and $s$. Applying eq. (\ref{eq. condensed}) to the particular solution $y^{par}_{t,s}$ in eq. (\ref{PARTICULAR SOLUION}), Theorem \ref{TheoGenSol} follows. 

\begin{theorem}\label{TheoGenSol} An equivalent explicit  representation of $y_{t}$ in eq. {\rm (\ref{TVAR(P)})} in terms of prescribed random variables $\{y_r\}_{s+1-p\le r\le s}$ for any  $s\in\integers$ and $t\in\integers_{s+1-p}$ such that  $s<t$ is given by: 
\begin{equation}  \label{TVAR(p)SOL1}
y_{t}=\underset{\text{Homogeneous Solution Part}}
{\underbrace{\displaystyle \sum_{m=1}^{p}\xi^{(m)}(t,s)y_{s+1-m}}}+\underset{\underset{\text{{\scriptsize Drift Part}}}{\text{Particular Solution:}}}{\underbrace{\sum_{r=s+1}^{t}\xi (t,r)\varphi (r)}+}\underset{\text{Particular Solution: Innovation Part}}{\underbrace{\sum \limits_{r=s+1}^{t}\xi_{q}(t,r)\varepsilon _{r}+\sum
\limits_{r=s+1-q}^{s}\xi _{s,q}(t,r)\varepsilon_{r}.}}
\end{equation}
\end{theorem}
The algebraic part of the proof, concerning TV-LDEs($p$) is provided in (Paraskevopoulos and Karanasos 2021), the stochastic part is given in the Appendix.
The right-hand side of eq. (\ref{TVAR(p)SOL1}) comprises four summation parts. In view of Proposition \ref{DEFKSIQ}, the sum of its last three parts is the particular solution given by eq. (\ref{PARTICULAR SOLUION}). More analytically, the first sum (the homogeneous solution in eq. (\ref{alternative homogeneous solution})) is a product of $m$ fundamental solutions multiplied by observable random variables. The second sum (the drift part of the particular solution in eq. (\ref{PARTICULAR SOLUION})) is formed by products involving the principal determinant $\xi (t,r)$  multiplied by the drift $\varphi(r)$. The terms of the third sum (the first part of the ``MA decomposition", see eq. (\ref{eq. condensed})) are the banded Hessenbergian coefficients $\xi _{q}(t,r)$,  multiplied by the unobservable errors. Finally, the terms of the fourth sum (the second part of the ``MA decomposition") are the banded Hessenbergian coefficients $\xi_{s,q}(t,r)$, multiplied by observable errors. In the Online Appendix G1, we apply eq. (\ref{TVAR(p)SOL1}) to an AR(1) model, recovering the well known explicit  representation of the latter.

\begin{remark}
\label{CORk=1} When $s=t-1$ (or $k=1$) the expression in Theorem {\rm \ref{TheoGenSol}} coincides with eq. {\rm (\ref{TVAR(P)})}. This is a consequence of the following statements: i) $\xi^{(m)}(t,t-1)=\phi_{m}(t)$ (see the discussion next to eq.  (\ref{Exdef:xi^m1}))
and ii) $\sum_{r=t}^{t}\xi(t,r)[\varphi(r)+u_{r}]=\varphi(t)+u_{t}$.  
\end{remark}

\begin{remark}
\label{Green's function representation} Replacing the homogeneous solution part in eq. {\rm (\ref{TVAR(p)SOL1})} by eq. {\rm (\ref{Unified homogeneous solution})},
and its innovation part with the left-hand side of eq. {\rm (\ref{eq. condensed})}, we can solely express eq. {\rm (\ref{TVAR(p)SOL1})} as a linear combination of the Green function. The expression in eq. {\rm (\ref{TVAR(p)SOL1})} has to be compared with the corresponding result in (Agarwal, 2000 p.77, eq. (2.11.8)), in which the fundamental solutions, denoted there by ``$v_{i}(k)$", are not in general explicitly expressed, but only in specific cases (see the examples in the previously cited reference).
\end{remark}

The methodology presented in this Section can be used in the study of infinite order autoregression models as well as in the case of the fourth
order moments for time-varying GARCH models. In the interest of brevity the detailed examination of the aforementioned models will be the subject of
future papers. We should also mention that another mathematical tool of constant use in difference equations is the generalized continuous fraction
approach (see, Van de Cruyssen, 1979). The concept of matrix continued fraction was introduced in Hallin (1984), whereas Hallin (1986) showed the
close connection between the convergence of matrix continued fractions and
the existence of dominated solutions for multivariate difference equations
of order two.

Another advantage of our TV-ARMA representation in eq. (\ref{TVAR(p)SOL1}), is its generality. That
is, in deriving it we do not make any assumptions on the time dependent coefficients. Therefore, it does not require a case by case treatment. In other words, we suppose that the law of evolution of the coefficients is unknown, in particular they may be stochastic (either stationary or non stationary) or deterministic. Therefore, no restrictions are imposed on the functional form of the time-varying autoregressive and moving average coefficients. In the non stochastic case the model allows for known abrupt changes, smooth changes and mixtures of them. If the changes are smooth the coefficients can depend on an exogenous variable $x_{t}$\ or $t$\ or both.
In the case of stochastically varying coefficients the model includes the generalized random coefficient (GRC) AR specification (see, for example, Glasserman and Yao, 1995, and Hwang and Basawa, 1998) as a special case or allows for Markov switching behaviour (see, for example, Hamilton, 1989 and 1994, chapter 22). In both aforementioned cases it allows for periodicity. We should also mention that the solution includes the case where the variable coefficients depend on the length of the series (see the example in
Section \ref{SUBSECEXAMPL}).

\subsection{Gegenbauer Functions as Hessenbergians}

We conclude this Section with an example. We show how the Gegengbauer
functions can be expressed as Hessenbergians. For a discussion of the
Gegenbauer processes and their applications to economics and finance see
Baillie (1996) and Dissanayake et al. (2018); see also, Giraitis and Leipus
(1995) and Caporale and Gil-Alana (2011).

\begin{example}
\label{ex. Gegenbauer} The Gegenbauer (or ultraspherical) functions, denoted
by $c_{j}^{(d)}(\phi)$ (hereafter, for notational simplicity we use $c_{j}$), are defined to be the coefficients in the power-series expansion of the
following equation: \vspace{-0.1in}
\begin{equation*}
(1-2\phi z+z^{2})^{-d}=\sum \limits_{j=0}^{\infty }c_{j}z^{j},
\end{equation*}
for $\left \vert z\right \vert \leq 1$, $\left \vert \phi \right \vert \leq 1 $, and $0<d<\frac{1}{2}$. 
\end{example}
The easiest way to compute $c_{j}$ (for $j\geq 2$) using computers is based on the solution of the following  second order homogeneous linear difference equation with coefficients functions of $j$
\begin{equation}\label{Gegenbauer}
c_{j}=2\phi \left( \frac{d-1}{j}+1\right) c_{j-1}-\left( 2\frac{d-1}{j}
+1\right) c_{j-2}\ \ \mathrm{for}\ j\geq 2,
\end{equation}%
subject to the initial condition values $c_{0}=1$ and $c_{1}=2\phi d$ (see, for example, Baillie, 1996, Chung, 1996, and the references therein).
The solution $c_{j}$  of eq. (\ref{Gegenbauer}) is given in the following Proposition. This follows from eq. (\ref{alternative homogeneous solution})  (or the homogeneous solution part of Theorem \ref{TheoGenSol}) applied with $t=j$,  $s=1$, $\phi_{1}(j)=2\phi( \frac{d-1}{j}+1), \phi_{2}(j)=-( 2\frac{d-1}{j}+1)$ (e.g. $\phi_{1}(2)=\phi (d+1), \phi_{1}(3)=2\phi\frac{d+2}{3}, \phi_{2}(3)=-\frac{2d+1}{3}$). Besides $y_{j,1}^{hom}$ is identified  with $c_j$. 
\begin{proposition}
The $j$-th (for $j\geq 2$) Gegenbauer coefficient is given by $c_{j}=\xi (j,1)c_{1}+\xi^{(2)} (j,1)c_{0}$ or equivalently (by replacing $c_{0}=1$ and $c_{1}=2\phi d$): 
\begin{equation}\label{Gegenbauer2}
  c_{j} =\xi (j,1)2\phi d+\xi^{(2)} (j,1),
\end{equation}
where $\xi (j,1)$ is a $(j-1)$-th order banded Hessenbergian (in this case tridiagonal matrix): 
\begin{equation*}
\xi (j,1)=\left \vert 
\begin{array}{cccccc}
\phi (d+1) & -1 &  &  &  &  \vspace{0.02in}\\ 
-\frac{2d+1}{3} & 2\phi \frac{d+2}{3} & -1 &  &  &  \vspace{0.07in}\\ 
& -\frac{d+1}{2} & \phi \frac{d+3}{2} & -1 &  &  \\ 
&  & \ddots & \ddots & \ddots &  \\ 
&  &  & -\left(2\frac{d-1}{j-1}+1\right) & 2\phi\left( \frac{d-1}{j-1}+1\right) & -1 \vspace{0.05in}\\ 
&  &  &  & -\left( 2\frac{d-1}{j}+1\right) & 2\phi \left( \frac{d-1}{j}+1\right)
\end{array}
\right \vert .
\end{equation*}
\end{proposition}
Recall that $\xi^{(2)} (j,1)$ differs from $\xi (j,1)$  only in the first two elements of the first column (see eq. (\ref{Def1: Phi^m_k})), which, in this case, are: $-d$ (or $\phi_{2}(2)$) and $0$, respectively. It follows from eq. (\ref{Exdef:xi^m1}) that $\xi^{(2)}(j,1)=-d\xi (j,2)$. Thus eq. (\ref{Gegenbauer2}) can be rewritten as: $c_j=\xi (j,1)2\phi d-d\xi (j,2)$. Next, we apply the multi-linearity property of determinants along the first column to express  $c_j$ in eq. (\ref{Gegenbauer2}) as a single Hessenbergian:

\begin{corollary}
The $j$-th Gegenbauer polynomial can be explicitly expressed as a $(j-1)$-th order
banded Hessenbergian: 
\begin{equation*}
c_{j}=\left \vert 
\begin{array}{cccccc}
d[2\phi^{2}(d+1)-1] & -1 &  &  &  &  \vspace{0.02in}\\ 
-2\phi d\frac{2d+1}{3} & 2\phi \frac{d+2}{3} & -1 &  &  & \vspace{0.05in}
\\ 
& -\frac{d+1}{2} & \phi  \frac{d+3}{2} & -1 &  &  \\ 
&  & \ddots & \ddots & \ddots &  \\ 
&  &  & -\left( 2\frac{d-1}{j-1}+1\right) & 2\phi \left( \frac{d-1}{j-1}%
+1\right) & -1 \\ 
&  &  &  & -\left( 2\frac{d-1}{j}+1\right) & 2\phi \left( \frac{d-1}{j}%
+1\right)%
\end{array}%
\right \vert .
\end{equation*}
\end{corollary}

\section{Asymptotic Stability}
\label{sec. Stability}
In this Section we apply our result in eq. (\ref{TVAR(p)SOL1}) to provide conditions ensuring the asymptotic stability\footnote{As pointed out by Grillenzoni (2000) stability is a useful feature of stochastic models because it is a sufficient (although not necessary) condition for optimal properties of parameter estimates and forecasts. Since model (\ref{TVAR(P)}) can be expressed in Markovian form, the stability condition allows many other stability properties, such as irreducibility, recurrence, regularity, non evanescence and tightness (see Grillenzoni, 2000 for details).} of a family of   TV-ARMA($p,q$) processes satisfying eq. (\ref{TVAR(P)}).

\subsection{Stability Conditions}
\label{StabilityCondition} 

The asymptotic stability problem is to provide sufficient conditions such that a class of stochastic processes solving eq. (\ref{TVAR(P)}) approaches a solution independently of the prescribed $p$ random variables (the effect of the prescribed random variables is gradually dying out) as $s\rightarrow -\infty$, that is when the homogeneous solution in eq. (\ref{TVAR(p)SOL1}) tends to zero, under a prescribed type of convergence. The explicit representation of the homogeneous solution in eq. (\ref{alternative homogeneous solution}) makes it possible to provide such type of conditions in Theorem \ref{theo: asympt. stability in L2} ensuring the $L_2$ convergence to zero of the homogeneous solution, that is
$y^{hom}_{t,s}\overset{L_2}{\to} 0$, as $s\rightarrow -\infty$, which means that $\lim_{s\to -\infty}||y^{hom}_{t,s}||_{L_2}=0$, or equivalently that $\lim_{s\rightarrow -\infty}\expectation\big(y^{hom}_{t,s}\big)^2=0$.
\begin{theorem}\label{theo: asympt. stability in L2} \textrm{i)}  
Let the autoregressive coefficients $\phi_{m}(t)$ be deterministic. If $\sup_t{|\phi_m(t)}|<\infty$ for each $m$ with $1\leq m\leq p$, then a sufficient condition for an $L_2$-bounded stochastic process $y_t$  (that is $\sup_t\expectation(y^2_t) <\infty$)  which solves eq. (\ref{TVAR(P)}) to be  asymptotically stable (in $L_2$ sense) is: $\lim_{s \rightarrow -\infty }\xi (t,s)=0$ for each $t$. \\
\textrm{ii)} Let the autoregressive coefficients $\phi_{m}(t)$ be stochastic. If $\sup_t\expectation(\phi _{m}^{2}(t))<\infty$ for each $m$ with $1\leq m\leq p$, then a sufficient condition for an $L_2$-bounded stochastic process $y_t$,  which solves eq. (\ref{TVAR(P)}), to be asymptotically stable (in $L_2$ sense)  is $\xi(t,s)\overset{L_2}{\rightarrow }0$, as $s\rightarrow -\infty$ for each $t$ (see, for
details, Proposition \ref{PROPABSSUM1} in Section 7.4).\footnote{Goldie and Maller (2000) offered sufficient conditions for the asymptotic stability of an AR($1$) model with stochastically varying coefficients, that is the a.s. convergence of the solution, i.e. $\sum\nolimits_{t=1}^{\infty} \phi _{1}(1)\phi _{1}(2)\cdots \phi _{1}(t-1)\varepsilon _{t}<\infty $, a.s. (see also Bougerol and Picard, 1992). Recently, in a multivariate setting, Erhardsson (2014) showed that the only sufficient condition is: $\left\vert \phi _{1}(1)\phi _{1}(2)\cdots \phi _{1}(t)\right\vert \overset{a.s.}{\rightarrow }0$ as $t\rightarrow \infty $, a result analogous to Theorem \ref{theo: asympt. stability in L2} (ii) (see also Section \ref{SECSTV} below).}
\end{theorem}
In the next Section, the asymptotic stability plays an essential role for  the existence and uniqueness of TV-ARMA($p,q$) solution processes in $L_2$.

The conditions in Theorem \ref{theo: asympt. stability in L2}(ii) include the ``bounded random walk" of Giraitis et al. (2014) also used by Petrova (2019). Properties such as stability characterize the statistical properties ($\sqrt{T}$ convergence and asymptotic normality, where $T$ is the sample size) of least squares (LS) and quasi-maximum likelihood (QML) estimators of the time-varying coefficients.\footnote{Azrak and M\'{e}lard (2006) have considered the asymptotic properties of quasi maximum likelihood estimators for a large class of ARMA\ models with time dependent coefficients and heteroscedastic innovations. The coefficients and the variance are assumed to be deterministic functions of time, and depend on a finite number of parameters which need to be estimated. Other researchers have also considered the statistical properties of maximum likelihood estimators for very general non stationary models. For example, Dahlhaus (1997) has obtained asymptotic results for a new class of locally stationary processes, which includes TV-ARMA processes (see Azrak and M\'{e}lard, 2006, and the references therein).}

In the time invariant case since $\xi (t,s)$ depends neither on $t$ nor on $s$ but only on their difference, that is the forecasting horizon, $k$, the stability condition in \ref{theo: asympt. stability in L2}(i) reduces to $\lim_{k\rightarrow \infty }\xi _{k}=0$, which holds if and only if all the roots $\lambda_{m}$ in eq. (\ref{Toeplitz}) lie inside the unit circle.

\subsection{Two Illustrative Examples}
\label{Two Illustrative Examples}
The first of the following examples concerns the logistic smooth transition
AR($1$) model (see, for example, Ter\"{a}svirta, 1994). The second example concerns periodic AR models.
\begin{example}
\label{logistic smooth transition} In the first model the autoregressive coefficient is given by : $\phi (t)=\phi_{1}F(t;\gamma ,\tau )+[1-F(t;\gamma ,\tau)]\phi_{2}$ (we drop the subscript $1$), where $F(t;\gamma ,\tau )=[1+e^{\gamma (t-\tau)}]^{-1}$, $\gamma \in \reals_{\geq 0}$, $\tau \in \integers$, is the first-order logistic function. Clearly, if $t>\tau $, then $F((t;\gamma ,\tau ))<0.5$ and regime $2$ prevails, whereas if $t<\tau$, then $F(t;\gamma ,\tau )>0.5$ and regime $1$ prevails. Let also $t_{2}$ be the value of $t$ for which $F(t;\gamma ,\tau )=0$ if $t\geq t_{2}$ and thus $\phi (t)=\phi_{2}$. Similarly, let $t_{1}$ be the value of $t$ for
which $F(t;\gamma ,\tau )=1$ if $t\leq t_{1}$ and thus $\phi (t)=\phi _{1}$
(clearly $t_{2}>t_{1}$, since $F(\cdot )$ is a decreasing function of time).
For this model, if $t\geq t_{2}$ and $s\leq t_{1}$ then $\xi (t,s)=\phi
_{1}^{t_{1}-s+1}\prod \nolimits_{r=t_{1}+1}^{t_{2}-1}\phi (r)\phi
_{2}^{t-t_{2}+1}$. Clearly, $\lim_{s\rightarrow -\infty }\xi (t,s)=0$ if and
only if $\left \vert \phi _{1}\right \vert <1$, whereas $\lim_{t\rightarrow
\infty }\xi (t,s)=0$ if and only if $\left \vert \phi _{2}\right \vert <1$.
\end{example}


\begin{example}
\label{ex. periodic AR1} In the periodic {\rm AR($1;\ell $)} model, $\ell \in \integers_{\geq 1}$ is the
number of seasons (i.e., quarters) and $\phi _{r}$, $r=1,\ldots, \ell $ denote the periodically varying autoregressive coefficients. Moreover, $T\in 
\integers_{\geq 0}$ stands for the number of periods (i.e., years). Accordingly $t=T\ell
+r$ is time expressed in terms of seasons (i.e., if $\ell =4$, $r=4$ and $T=1$, then $t=8$ quarters). If we want to forecast $k\ell $ seasons ahead,
that is $t-s=k\ell $ or $s=t-k\ell $, then: $\xi (t,s)=[\prod \nolimits_{r=1}^{\ell }(\phi _{r})]^{k}$. Clearly $\left \vert \phi
_{r}\right \vert <1$ for all $r $ is a sufficient but not necessary condition for $\lim_{s\rightarrow -\infty }\xi (t,s)=0$ (or equivalently $\lim_{k\rightarrow \infty }\xi (t,s)=0$). The necessary and sufficient condition is $\left \vert \prod \nolimits_{r=1}^{\ell }(\phi_{r})\right
\vert <1$.\footnote{For a study of the periodic stationarity of a random coefficient periodic
autoregression (RCPAR) see, for example, Aknouche and Guerbyenne (2009).}
\end{example}

\section{Solution Representation in $L_2$}
\label{Representations in L_2}

Having specified a general method for manipulating explicitly the TV-ARMA type of models, we turn our attention to a consideration of their Wold-Cram\'{e}r decomposition followed by the fundamental second order properties of these models. In this Section we shall restrict ourselves to the treatment of DTV-ARMA processes (those with deterministic coefficients). In Section \ref{SECSTV} we present STV-AR processes (those with stochastically varying coefficients), which incorporate the GRC and double stochastic AR models.

\subsection{Wold-Cram\'{e}r Decomposition}
\label{WoldCramerDecomposition}
In Theorem \ref{WoldCramDec}, we provide the existence of the Wold-Cram\'{e}r decomposition  (see Cram\'{e}r, 1961)\footnote{Since a non-stationary generalization of Wold's result was given by Cram\'{e}r, it is referred to as Wold-Cram\'{e}r decomposition.} and, therefore, impulse response functions (IRFs), for the DTV version of the model in eq. (\ref{TVAR(P)}). In particular, we provide an explicit  condition  (the absolute summability of the principal determinant): 
\begin{equation}\label{abs summability condition}
\sum_{r=-\infty }^{t}|\xi (t,r)|<\infty ,\ \text{for\ all}\ t,\ \ \ \ \ (\text{absolute\ summability\ condition})  
\end{equation}%
which, along with the boundedness of the drift and the moving average coefficients, ensure the existence of the Wold-Cram\'{e}r decomposition (Theorem \ref{WoldCramDec}) of DTV-ARMA models with finite first two unconditional moments and autocovariance function (Propositions \ref{ABSSUM}, \ref{ProSecMom}). 
In Theorem \ref{uniqueness of the Wold-Cramer decomposition} it is shown that every asymptotically stable stochastic process, which solves eq. (\ref{TVAR(P)}) coincides with the Wold-Cram\'{e}r decomposition in eq. (\ref{Wold Cramer}). 


As the absolute summability condition in (\ref{abs summability condition}) also guarantees the absolute summability of $\xi_{q}(t,r)$ in eq. (\ref{KSIq}) the following Theorem holds:

\begin{theorem}[Existence]\label{WoldCramDec}
Let the absolute summability condition in  {\rm (\ref{abs summability condition})} hold. Let also $\varphi(t)$ and $\theta_{l}(t)$ be  bounded functions in $t$. Then there exists a solution of eq. {\rm (\ref{TVAR(P)})} in $L_2$ of the form
\begin{equation}  \label{Wold Cramer}
y_{t}=\sum_{r=-\infty }^{t}\xi (t,r)\varphi (r)+\displaystyle\sum_{r=-\infty }^{t}\xi _{q}(t,r)\varepsilon_{r},
\end{equation}
that is an explicit representation of the Wold-Cram\'{e}r solution decomposition.
\end{theorem}
A direct proof of Theorem \ref{WoldCramDec} is given in the Appendix ensuring that $y_t$ in eq. (\ref{Wold Cramer}) is second order, that is $y_t\in L_2$.  In particular, show there that $\sup_t\expectation(\upsilon^2_t)<\infty$ (for the definition of $\upsilon_{t}$ see  eq.  (\ref{Forcing term})), which extends a result established by Peiris (1986) (see the proof of his Theorem 2.1, provided that $\varphi(t)=0$). In addition, we deduce the equivalence between $\sum_{r=-\infty }^{t}\xi^2_{q}(t,r)<\infty$ and $y_t\in L_2$, provided that $\sum_{r=-\infty }^{t}\xi (t,r)\varphi (r)$ converges in $\reals$ (see Proposition G1, in the Online Appendix G2). This  extends a result mentioned in the above cited reference, that is $y_t$ is second order if and only if eq. ``(2.4)" holds true, provided that $\varphi(t)=0$.
The solution of eq. (\ref{TVAR(P)}) in eq. (\ref{Wold Cramer}) is decomposed into two orthogonal parts (see for a proof the Additional Appendix L1), a deterministic part and a mean zero random part, that is, $\expectation(y_{t})=\sum_{r=-\infty }^{t}\xi (t,r)\varphi (r)$ is the non random part (see Proposition \ref{ABSSUM} below), while $\lim_{s\to -\infty} \mathbb{FE}_{t,s}=\sum_{r=-\infty }^{t}\xi _{q}(t,r)\varepsilon_{r}$ (i.e., the limit of forecast errors, see eqs. (\ref{FE and MSE1}) below) is the mean zero random part of $y_{t}$. 
As eq. (\ref{Wold Cramer}) is future independent, we shall also referred to it as a causal solution of DTV-ARMA models.

Hallin (1979), Singh and Peiris (1987), Kowalski and Szynal (1991), Grillenzoni (2000), and Azrak and M\'{e}lard (2006) evaluate the Wold-Cram\'{e}r decomposition through recursion. In sharp contrast, eq. (\ref{Wold Cramer}) in Theorem \ref{WoldCramDec} provides a direct and explicit representation of the Wold-Cram\'{e}r decomposition expressed in terms of banded Hessenbergians with entries the coefficients of the model in eq. (\ref{TVAR(P)}). Moreover, in the Appendix we show that the Wold-Cram\'{e}r decomposition in eq. (\ref{Wold Cramer}) can be expressed in a  compact form:  $$y_t=\sum_{r=-\infty }^{t}\xi (t,r)\upsilon_r.$$ 

Let $\varphi(r)=0$ and $\{\varepsilon_{t}\}_t$ be an orthonormal sequence in $L_2$, that is  $\expectation(\varepsilon_{t}\varepsilon_{s})=\delta_{t,s}$ (Kroneker delta). Under these additional assumptions to those considered in Theorem \ref{WoldCramDec}, the latter recovers the following uniqueness result: Eq. (\ref{Wold Cramer}) is the unique solution  of eq. (\ref{TVAR(P)}) (see Hallin and Ingenbleek, 1983, and Peiris, 1986). In the Online Appendix G1 we apply eq. (\ref{Wold Cramer}) to an AR(1) model, recovering the well known Wold solution representation of the latter.

In the Appendix we show that eq. (\ref{Wold Cramer}) can be alternatively derived from  eq. (\ref{TVAR(p)SOL1}), when the latter is applied with prescribed random variables
 $$y_{s+1-m}=\sum_{r=-\infty }^{s+1-m}\xi (s+1-m,r)\upsilon_{r}\ \ \ \text{for}\ \ \ 1\leq m\leq p,$$ 
for any arbitrary $s$ such that $s<t$. This result arises from the formulas
\begin{equation}\label{L2 convergence of hom. solution}
\left. \begin{array}{l} 
\displaystyle\lim_{s\rightarrow -\infty }\sum_{m=1}^{p}\xi ^{(m)}(t,s)\displaystyle\sum_{r=-\infty }^{s+1-m}\xi (s+1-m,r)\upsilon_{r}\overset{L_{2}}{=}0,\\\\ \displaystyle\lim_{s\rightarrow -\infty }\displaystyle\sum_{r=s+1-q}^{s}\xi
_{s,q}(t,r)\varepsilon _{r}\overset{L_{2}}{=}0
\end{array}\right\}  
\end{equation}
as follows: 
Applying  eqs. (\ref{L2 convergence of hom. solution}) to eq. (\ref{TVAR(p)SOL1}), it follows from eq. (\ref{PARTICULAR SOLUION}) that  eq. (\ref{Wold Cramer}) can be equivalently rewritten as: $y_{t}=\lim_{s\rightarrow -\infty}y_{t,s}^{par}$. The first of the formulas in (\ref{L2 convergence of hom. 
solution}) ensures that the Wold-Cram\'{e}r solution decomposition in eq. (\ref{Wold Cramer}) is asymptotically stable. Moreover, in the Appendix we show that every  asymptotically stable stochastic process in $L_2$, which solves eq. {\rm (\ref{TVAR(P)})} coincides almost surely with $y_t$ in eq. {\rm (\ref{Wold Cramer})}, as demonstrated in the following Theorem: 
\begin{theorem}[Uniqueness]\label{uniqueness of the Wold-Cramer decomposition}
Let the conditions of Theorem \ref{WoldCramDec} hold. Then there exists a uniquely determined  asymptotically stable stochastic process in $L_2$, which solves eq. {\rm (\ref{TVAR(P)})} that  coincides almost surely with $y_t$ in eq. {\rm (\ref{Wold Cramer})}.
\end{theorem}
Assuming in addition to the conditions of Theorem \ref{uniqueness of the Wold-Cramer decomposition} that $\sup_t{|\phi_m(t)|}<\infty$ for each $m$ with $1\leq m\leq p$ and taking into account that the absolute summability condition in (\ref{abs summability condition}) implies that $\lim_{s\to -\infty}\xi(t,s)= 0$, we infer from Theorem \ref{theo: asympt. stability in L2}(i) that every bounded process that solves eq. (\ref{TVAR(P)}), as being asymptotically stable, also coincides with the Wold-Cram\'{e}r decomposition in eq. (\ref{Wold Cramer}). The existence of such a bounded process is guaranteed by the boundedness of the absolute summability, as a function of $t$, that is
\begin{equation}\label{abs summability condition2}
\sum_{r=-\infty }^{t}|\xi (t,r)|\le C <\infty,\ \text{for\ all}\ t\ \ \text{and\ some}\ \  C\in \reals_{> 1}.     
\end{equation}
We show in the following Corollary that under the condition (\ref{abs summability condition2}), the Wold-Cram\'{e}r solution decomposition in eq. (\ref{Wold Cramer}) is bounded. Moreover, as (\ref{abs summability condition2}) guarantees the conditions of Theorem \ref{theo: asympt. stability in L2}(i) (see for a proof the Proposition L1 in the Additional Appendix L2), it arises from Theorem \ref{uniqueness of the Wold-Cramer decomposition} that the following uniqueness result holds: 
\begin{corollary}\label{uniqueness of the Wold-Cramer decomposition2}
Let the conditions of Theorem \ref{WoldCramDec} hold. Let also {\rm (\ref{abs summability condition2})} holds. Then there exists a uniquely determined $L_2$-bounded stochastic process, which solves eq. {\rm (\ref{TVAR(P)})} and coincides almost surely with $y_t$ in eq. {\rm (\ref{Wold Cramer})}. 
\end{corollary}
A proof of these results is included in Appendix. Corollary \ref{uniqueness of the Wold-Cramer decomposition2} must be compared with an analogous result  established by Neimi (Theorem 2.1, 1983) for ARMA  non-stationary processes with constant coefficients coupled with the AR-regularity condition and zero drift ($\varphi(t)=0$). An extension of this result covering non-stationary TV-ARMA($p,q$) models coupled with zero drift and extended regularity conditions was obtained by Kowalski and Szynal (1988).

Kowalski and Szynal (1991) and Grillenzoni (2000) derived sufficient conditions ensuring that $y_t$ is  second order, that is $y_t\in L_2$ (Kowalski and Szynal examined the case of deterministic coefficients with zero drift, whereas Grillenzoni allowed the coefficients to be stochastic as well, see also Remark 5 in Section 7.4). 
These conditions for a TV-AR($p$) model with zero drift are summarized in the following Remark.

\begin{remark}
\label{KS&G} The above mentioned sufficient conditions for $y_t\in L_2$ are as follows: i) The deterministically varying polynomial $\Phi
_{t}(z^{-1})=1-\sum_{m=1}^{p}\phi _{m}(t)z^{-m}$ is regular. That is, $\phi	_{m}(t)$ are such that there exist the limits $\lim_{t\rightarrow \infty
}\phi _{m}(t)=\phi _{m}$ and $\sum\nolimits_{r=1}^{\infty }\varrho 	^{2r}<\infty $, where $\varrho =\varrho (\Phi )+\epsilon $, $\epsilon >0$, $\varrho (\Phi )=\max \{\left\vert z_{m}\right\vert ,$ $\Phi (z_{m}^{-1})=0\}$ with $\Phi (z^{-1})=1-\sum_{m=1}^{p}\phi _{m}z^{-m}$ (see eq. (8) in Kowalski and Szynal, 1991).\footnote{Kowalski and Szynal (1991) showed that $\varrho (\Phi )$ is the spectral
radius of the matrix $\mathbf{\Phi =}\left[ 
\begin{array}{ccccc}
		0 & 0 & \cdots  & 0 & -\phi _{m} \\ 
		1 & 0 & \cdots  & 0 & -\phi _{m-1} \\ 
		0 & 1 & \cdots  & 0 & -\phi _{m-2} \\ 
		\vdots  & \vdots  & \vdots \vdots \vdots  & \vdots  & \vdots  \\ 
		0 & 0 & \cdots  & 1 & -\phi _{1}\vspace{-0.1in}%
		\end{array}
		\right] $\vspace{0.15in}\\

		(see page 75 in their paper).} ii) The deterministically varying polynomial $\Phi _{t}(z^{-1})$ should have roots that entirely lie inside the unit circle, with the exception, at most, of a finite set of points (see Proposition 1 in Grillenzoni, 2000).
\end{remark}

The sufficient conditions in Remark \ref{KS&G} are not, however, necessary, since they do not cover the case of periodic coefficients, see Grillenzoni (1990) or Karanasos et al. (2014,a,b). Whereas in the case of TV-AR($p$) models  $\xi_q(t,r)=\xi(t,r)$, we infer from
the discussion next to Theorem \ref{WoldCramDec} (see also Proposition G.2 in the Online Appendix G) the equivalence between $\sum_{r=-\infty }^{t}\xi^2(t,r)<\infty$ (square summability) and $y_t\in L_2$, provided that $\varphi(t)=0$. As a consequence, the stability condition $\lim_{s\rightarrow -\infty }\xi (t,s)=0$ is necessary for the square summability to hold or equivalently for $y_t\in L_2$  and therefore it is also necessary for the two conditions in Remark \ref{KS&G}, as demonstrated below:
$$ \{\text{The Sufficient Conditions in Remark \ref{KS&G}}\} \implies y_t\in L_2 \iff  \sum_{r=-\infty }^{t}\xi^2(t,r)<\infty \implies \lim_{s\rightarrow -\infty }\xi (t,s)=0. $$



\subsection{Unconditional Moments}\label{UNCONDMOMENTS}
In this Subsection we present explicit formulae for the first and second unconditional moments of the Wold-Cram\'{e}r solution decomposition of the DTV-ARMA family of processes coupled with sufficient and necessary conditions for their existence, as demonstrates the following Proposition.

\begin{proposition}
\label{ABSSUM} Let the conditions of Theorem \ref{WoldCramDec} hold. Then the unconditional mean of the process $y_{t}$ in eq. {\rm (\ref{Wold Cramer})} with deterministic coefficients, exists in $\reals$ and is given by 
\begin{equation}\label{first unconditional moment}
\expectation(y_{t})=\sum\limits_{r=-\infty }^{t}\xi (t,r)\varphi (r).
\end{equation}
The absolute summability condition is also sufficient for the existence of the unconditional variance of the above mentioned process $y_{t}$, which is given by 
\begin{equation}\label{VAR}
\var(y_{t})=\sum_{r=-\infty }^{t}\xi _{q}^{2}(t,r)\mathbb{\sigma}^{2}(r).  
\end{equation}
Necessary conditions for the  $y_t$ process to be  second order respectively are: 
\begin{equation*}
\lim_{s\rightarrow -\infty }\xi (t,s)\varphi (s)=0\ \mathrm{and}\
\lim_{s\rightarrow -\infty }\xi _{q}^{2}(t,s)\sigma ^{2}(s)=0\ \mathrm{for\
all}\ t.
\end{equation*}%
Moreover, the stability condition, that is $\lim_{s\rightarrow -\infty }\xi(t,s)=0$, is sufficient for the above two limits to exist, due to the boundedness of $\varphi (r)$ and $\sigma ^{2}(r)$, while it is necessary for the absolute summability to hold.
\end{proposition}

A proof of Proposition \ref{ABSSUM} is provided in Appendix, noticing that the unconditional mean $\expectation(y_{t})$ is the same for both the AR and the ARMA processes.

The main logical connections between the conditions, described in the above Proposition, are summarized in the following commutative diagrams:

\begin{figure}[th]
\begin{equation}\label{commutative diagram}
\begin{array}{ccccc}
\displaystyle\sum\limits_{r=-\infty }^{t}\xi (t,r)\varphi (r)\in \reals
& \Longleftarrow & \displaystyle\sum_{r=-\infty }^{t}|\xi (t,r)|<\infty & 
\Longrightarrow & \displaystyle\sum_{r=-\infty }^{t}\xi _{q}^{2}(t,r)\mathbb{%
\sigma }^{2}(r)<\infty \vspace{-0.1in} \\ 
\rotatebox{-90}{$\Longrightarrow$} & \rotatebox{-135}{$\Longrightarrow$} & %
\rotatebox{-90}{$\Longrightarrow$} & \rotatebox{-45}{$\Longrightarrow$} & %
\rotatebox{-90}{$\Longrightarrow$}\vspace{0.1in} \\ 
\lim\limits_{s\rightarrow -\infty }\xi (t,s)\varphi (s)=0 & \Longleftarrow & 
\lim\limits_{s\rightarrow -\infty }\xi (t,s)=0 & \Longrightarrow & 
\lim\limits_{s\rightarrow -\infty }\xi _{q}^{2}(t,s)\sigma ^{2}(s)=0.%
\end{array}
\end{equation}%
\vspace{-0.12in}
\caption{Commutative Diagrams}
\end{figure}

\subsection{Autocovariance Function}
\label{Autocovariance Function}
In the following Proposition, we state an explicit expression for the covariance structure for the  Wold-Cram\'{e}r solution decomposition of the DTV-ARMA($p,q$) process.

\begin{proposition}
\label{ProSecMom} Let the conditions of Theorem \ref{WoldCramDec} hold. Then the time
varying $\ell$-order autocovariance function $\gamma_{t}(\ell)=\autocov(y_{t},y_{t-\ell})$, $\ell \in \integers_{\geq0}$, of $y_t$ in eq. {\rm (\ref{Wold Cramer})}, exists in $\reals$ and is given by 
\begin{equation}  \label{autocovariance function}
\gamma_{t}(\ell)=\sum_{r=-\infty}^{t-\ell}\xi_{q}(t,r)\xi_{q}(t-\ell ,r)\mathbb{\sigma}^{2}(r).
\end{equation}
\end{proposition}%
A proof of Proposition \ref{ProSecMom} is found in the Appendix. The time-varying variance of $y_{t}$ in eq. (\ref{VAR}), is recovered by
applying $\gamma_{t}(\ell)$ for $\ell=0$, that is $\gamma_{t}(0)=\var(y_{t})$. Moreover, the absolute summability condition implies absolute
summable autocovariances: $\sum \nolimits_{\ell=0}^{\infty}\left
\vert \gamma_{t}(\ell)\right \vert <\infty$ for all $t$. Notice again that for a pure AR process, the autocovariance function formula in eq. (\ref{autocovariance function}) must be applied for $\xi_{0}(t,r)=\xi(t,r)$.

From a computational viewpoint, the covariance structure of $\{y_{t}\}_t$ can be numerically evaluated by computing the banded Hessenbergian coefficients, $\xi_{q}(t,r)$ in eq. (\ref{KSIq}) and substituting these in eq. (\ref{autocovariance function}).

The next remark highlights the importance of the existence of finite second
moments.

\begin{remark}
Azrak and M\'{e}lard (2006) considered the asymptotic properties of QML
estimators for the DTV-ARMA family of models where the coefficients depend
not only on $t$ but on $T$ as well (see Alj et al., 2017, for the
multivariate case). In their Theorem and Lemma $1$ the existence of finite
second moments was required. They also show that the dependence of the model
with respect to $T$ has no substantial effect on their conclusions except
that a.s. convergence is replaced by convergence in probability, since convergence in $L_{2}$ norm implies convergence in probability (see Lemma $1^{\prime }$ in their paper).
\end{remark}

We conclude this Section with two more examples and a discussion of forward asymptotic stability.

\paragraph{Two More Examples\label{SUBSECEXAMPL}:}

The following examples, concerning AR($1$) processes with variable autoregressive coefficients, are taken from Azrak and M\'{e}lard (2006).

\begin{example}
\label{ex. periodic autoregressive1} In the first example, the autoregressive coefficient is a periodic function of time defined by 
\begin{equation*}
y_{t}=\phi (t)y_{t-1}+\varepsilon _{t},
\end{equation*}%
where $\varepsilon _{t}$ is a martingale difference sequence defined on $L_{2}$ with
constant variance $\sigma ^{2}$. Moreover, the autoregressive coefficient is
given by $\phi (t)=\beta _{t-n[t/n]}$, where $n\in \integers_{\geq 2}$ and $%
[x]$ is the larger integer less or equal to $x$ (see also Dahlhaus, 1996).
By specializing the results of Theorem {\rm \ref{WoldCramDec}} and Proposition {\rm \ref{ABSSUM}}, the Wold-Cr\'{a}%
mer decomposition (if and only if $\left\vert \beta \right\vert <1$, where $%
\beta =\beta _{0}\cdot \beta _{1}\cdots \cdot \beta _{n-1}$) is given by 
\begin{equation*}
y_{t}=\sum\limits_{r=-\infty }^{t}\xi (t,r)\varepsilon _{r},
\end{equation*}%
with unconditional variance 
\begin{equation*}
\var(y_{t})=\sigma ^{2}\sum\limits_{r=-\infty }^{t}\xi ^{2}(t,r),
\end{equation*}%
where 
\begin{equation*}
\xi (t,r)=\beta ^{\lbrack \frac{t-r}{n}]}\left( \prod\limits_{j=0}^{t-r-1-n[%
\frac{t-r}{n}]}\beta _{t-j}-n\left[ \frac{t-j}{n}\right] \right) ,
\end{equation*}%
and, therefore 
\begin{equation*}
\sum\limits_{r=-\infty }^{t}\xi ^{2}(t,r)=\frac{1}{1-\beta ^{2}}%
\sum\limits_{r=t-n+1}^{t}\left( \prod\limits_{j=0}^{t-r-1-n[\frac{t-r}{n}%
]}\beta _{t-j}-n\left[ \frac{t-j}{n}\right] \right)
\end{equation*}%
(see also eq. (4.2) in Azrak and M\'{e}lard, 2006).
\end{example}

\begin{example}
In the second example (see Example 2 in Azrak and M\'{e}lard, 2006), the
autoregressive coefficient is an exponential function of time given by 
\begin{equation*}
\phi (t)=\left \{ 
\begin{array}{lll}
\phi & \text{for} & t\leq 0, \\ 
\phi \lambda ^{t/T} & \text{for} & t=1,\ldots ,T-1, \\ 
\phi \lambda & \text{for} & t\geq T,%
\end{array}
\right.
\end{equation*}
where $T\in \integers_{\geq 1}$ is the sample size. For this case (assuming
that $t>T$) 
\begin{equation*}
\xi (t,r)=\left \{ 
\begin{array}{lll}
(\phi \lambda )^{t-r} & \text{for} & r\in \lbrack T,t], \\ 
(\phi \lambda )^{t}\phi ^{-r}\lambda ^{-\left( \frac{T+1}{2}+\frac{r(r-1)}{%
2T }\right) } & \text{for} & r=1,\ldots ,T-1, \\ 
\phi ^{1-r}\xi (t,1) & \text{for} & r\leq 0.%
\end{array}
\right.
\end{equation*}
The condition $\left \vert \phi \right \vert <1$ (necessary and sufficient)
entails: 
\begin{equation*}
\sum \limits_{r=-\infty }^{t}\xi ^{2}(t,r)=\frac{1}{1-\phi ^{2}}\xi
^{2}(t,1)+(\phi \lambda )^{2t}\sum \limits_{r=1}^{T-1}\phi ^{-2r}\lambda
^{-(T+1+\frac{r(r-1)}{T})}+\frac{1-(\phi \lambda )^{2(t-T+1)}}{1-(\phi
\lambda )^{2}}.
\end{equation*}
As pointed out by Azrak and M\'{e}lard (2006) the use of variable
coefficients, which depend on the length of the series, is compatible with
the approach of Dahlhaus (1997).
\end{example}

\section{Forecasting and Asymptotic Efficiency}\label{SecSecMom}
In this Section we provide explicitly a sufficient condition for the invertibility of a causal  DTV-ARMA$(p,q)$ process. In conjunction with Theorem \ref{WoldCramDec}, the above mentioned condition makes it possible to obtain in Subsection \ref{SUBPREDICTOR1} an explicit form to the $k$-step-ahead optimal (in $L_{2}$ sense) linear predictor along with the associated MSE (mean square error). An alternative optimal forecast is based on a finite sequence of observable random variables, employing for that the representation given by eq. (\ref{TVAR(p)SOL1}). Both approaches yield identical forecast and mean square errors. In this context, in Subsection \ref{FORWASYMPTEFFIC} we discuss, in a unified scheme,  the forward asymptotic efficiency of  the model with deterministic coefficients  given by eq. (\ref{TVAR(P)}).

\subsection{Invertibility\label{Invertibility}}
A DTV-ARMA$(p,q)$ process is invertible if and only if the current value of the input $\varepsilon _{t}$ can be expressed as a (converging) linear combination of
the present and past random variables $y_r$ ($r\le t$) (see Brockwell and Davis, 2016, p.76). The main result of this Section is presented in Theorem \ref{Invertibility1} (see below).

Eq. (\ref{TVAR(P)}) can be rewritten as\vspace{-0.1in} 
\begin{equation}
\varepsilon_{t}=y_{t}-\varphi (t)-\sum \limits_{l=1}^{q}\theta _{l}(t)\varepsilon
_{t-l}-\sum \limits_{m=1}^{p}\phi _{m}(t)y_{t-m}.  \label{EPSILON}
\end{equation}%
The principal matrix associated with the moving average part is defined by

\begin{equation}
\mathbf{\Theta }(t,s)\!=\!\!\left[ 
\!\!\begin{array}{cccccccc}
-{\footnotesize \theta }_{1}{\footnotesize (s+1)} & {\footnotesize -1} &  & 
&  &  &  &  \\ 
-{\footnotesize \theta }_{2}{\footnotesize (s+2)} & -{\footnotesize \theta }%
_{1}{\footnotesize (s+2)} & \ddots &  &  &  &  &  \\ 
{\footnotesize \vdots } & {\footnotesize \vdots } & \ddots & \ddots &  &  & 
&  \\ 
-{\footnotesize \theta }_{q}{\footnotesize (s+q)} & -{\footnotesize \theta }%
_{q-1}{\footnotesize (s+q)} & \ddots & \ddots & \ddots &  &  &  \\ 
& -{\footnotesize \theta }_{q}{\footnotesize (s+q+1)} & \ddots & \ddots & 
\ddots & \ddots &  &  \\ 
&  & {\small \ddots } & {\small \ddots } & {\small \ddots } & {\small \ddots 
} & \ddots &  \\ 
&  &  & -{\footnotesize \theta }_{q}{\footnotesize (t-1)} & {\footnotesize %
\theta }_{q-1}{\footnotesize (t-1)} & {\small \cdots } & -{\footnotesize %
\theta }_{1}{\footnotesize (t-1)} & {\footnotesize -1} \\ 
&  &  &  & -{\footnotesize \theta }_{q}{\footnotesize (t)} & {\small \cdots }
& -{\footnotesize \theta }_{2}{\footnotesize (t)} & -{\footnotesize \theta }%
_{1}{\footnotesize (t)}%
\end{array}\!\!
\right]  \label{Theta}
\end{equation}%
(for $k=t-s\geq 1$). The matrix $\mathbf{\Theta }(t,s)$ has a similar structure to the principal matrix associated with the AR operator, $\mathbf{\Phi }(t,s)$, that is both matrices are banded lower Hessenberg of order $k$. It
is clear that for $q\geq t-s$, $\boldsymbol{\Theta }(t,s)$ is a full lower Hessenberg matrix.

For every pair $(t,s)\in \integers^{2}$ with $t-s\geq1$ we define the principal determinant associated with eq. (\ref{Theta}): 
\begin{equation*}
\vartheta(t,s)=\det(\mathbf{\Theta}(t,s)).
\end{equation*}
Formally $\vartheta(t,s)$ (and similarly to $\xi(t,s)$), is a banded Hessenbergian. We further extend the definition of $\vartheta(t,s)$ by assigning the initial values: $\vartheta(s,s)=1$ and $\vartheta(t,s)=0$ for $t<s$. Accordingly, $\vartheta(t,s)$ is the Green function associated with the MA($q$) operator.

In analogy with the definition of $\xi_{q}(t,r)$ (see eq. (\ref{KSIq})), we
define: 
\begin{equation*}
\vartheta _{p}(t,r)=\vartheta (t,r)-\sum \limits_{m=1}^{p}\vartheta
(t,r+m)\phi _{m}(r+m).
\end{equation*}

In the following Theorem we give a sufficient condition for a causal DTV-ARMA($p,q$) process determined by eq. (\ref{Wold Cramer}) to be invertible. 

\begin{theorem}
\label{Invertibility1}  Let the conditions in Theorem \ref{WoldCramDec} hold and $\phi_{m}(t)$ be bounded in $t$. Let additionally the following absolute summability condition $\sum_{r=-\infty}^{t}|\vartheta (t,r)|<\infty$ hold for each $t$. Then the process in eq. (\ref{Wold Cramer}) is invertible, that is  
\begin{equation}  \label{invertibility}
\varepsilon_{t}=\displaystyle\sum_{r=-\infty }^{t}\vartheta _{p}(t,r)y_{r}-\sum_{r=-\infty }^{t}\vartheta (t,r)\varphi (r)
\end{equation}
and solves eq. {\rm (\ref{EPSILON})}. 
\end{theorem}
The proof of Theorem \ref{Invertibility1} essentially repeats the arguments of the  proof of Theorem \ref{WoldCramDec}; switching the
roles of $y_{r}$ and $\varepsilon_{r}$, and replacing $\xi (t,r)$ with $\vartheta (t,r)$, $\xi _{q}(t,r)$ with $\vartheta_{p}(t,r)$ and $\varphi(r)$ with $-\varphi(r)$. This result is fulfilled by showing in the Online Appendix H, that both terms in the right-hand side of eq. (\ref{invertibility}) are processes in $L_2$. Eq. (\ref{invertibility}) recovers the formula devised by Hallin (1979) in Theorem 3, but this time in a more explicit form.

Following laborious research work, the literature contains a diversity of ``time-varying" specifications of linear form whose main time series properties either remain unexplored or have not been fully examined. Making progress in interpreting seemingly different models we have provided in this Section a common platform for the investigation of their time series properties. 
With the help of a few detailed examples, i.e., smooth transition AR processes, periodic and cyclical
formulations, we have demonstrated how to encompass various time series processes within our unified framework. The significance of our methodology is almost self-evident from the large number of problems that it can solve. Our proposed approach allows us to handle ``time-varying" models of infinite order, by introducing unbounded order linear difference equation of index $p$. This type of equations also yield Hessenbergian solutions, but involving full lower Hessenberg matrices. The latter enables us to obtain easily handled explicit solutions for the infinite order case along with the fundamental properties of corresponding models. An advantage of our technique is that it can be applied with ease, that is without any major alterations, in a multivariate setting and provides a solution to the problem at hand without adding complexity. 

\subsection{Optimal Forecasts}
\label{SUBPREDICTOR1} 
Theorems \ref{WoldCramDec} and \ref{uniqueness of the Wold-Cramer decomposition} ensure the existence of a unique asymptotically stable solution $y_{t}$ of eq. (\ref{TVAR(P)}) given by eq. (\ref{Wold Cramer}). Let $k=t-s$ for $s<t$ be the forecasting horizon. Let also $\mathcal{M}_{s}$ be the smallest closed subspace of $L_{2}$ based on the sequence of past observations $\{y_{s},y_{s-1},...\}$, which also contains all constant functions (see Brockwell and Davis, 1991, p. 62-64).  Assuming further the invertibility conditions in Theorem  \ref{Invertibility1}, as the terms of the sequence of errors $\{\varepsilon_{s},\varepsilon_{s-1}\ldots\}$ can be expressed in terms of $y$'s by eq. (\ref{invertibility}), they also belong to the space $\mathcal{M}_{s}$. The optimal predictor of $y_{t}$ in eq. (\ref{Wold Cramer}), i.e., its  conditional expectation on $\mathcal{M}_{s}$, coincides with its linear projection on the space spanned by $\{y_{s},y_{s-1},...\}\cup\{1\}$, whenever  $\{\varepsilon_r\}_r$ is a martingale difference sequence (see Brockwell and Davis, 2016, p. 334, see also Rosenblatt, 2000, p. 83).

The conditional expectation of $y_{t}$ given $\mathcal{M}_{s}$, that is $\expectation(y_{t}|\mathcal{M}_{s})$, turns out to
be the orthogonal projection of $y_{t}$ on $\mathcal{M}_{s}$, an explicit form of which is given in the following Proposition.

\begin{proposition}
\label{Forecasting infinite observations} The $k$-step-ahead optimal (in $L_{2}$-sense) predictor of  $y_t$, as stated above, is linear and it is given by 
\begin{equation}  \label{optimal linear predictor}
\expectation(y_{t}|\mathcal{M}_{s})=\sum_{r=-\infty }^{t}\xi (t,r)\varphi (r)+\displaystyle \sum_{r=-\infty}^{s}\xi_{q}(t,r)\varepsilon_{r},
\end{equation}
where $\varepsilon _{r}$ for $-\infty< r\leq s$ are given by eq. {\rm (\ref{invertibility})}.  The forecast error associated with $\expectation(y_{t}|\mathcal{M}_{s})$, that is $\mathbb{FE}_{t,s}\overset{\mathrm{def}}{=}y_{t}-\expectation (y_{t}\left\vert
{}\right. \mathcal{M}_{s})$, and the associated mean square error $\mathbb{M}\mathbb{SE}_{t,s}\overset{\mathrm{def}}{=}\var(\mathbb{FE}_{t,s})$ are  given by: \vspace{-0.12in} 
\begin{equation}  \label{FE and MSE1}
\mathbb{FE}_{t,s}=\sum_{r=s+1}^{t}\xi _{q}(t,r)\varepsilon _{r},\ \mathbb{MSE}_{t,s}=\sum_{r=s+1}^{t}\xi _{q}^{2}(t,r)\mathbb{\sigma }^{2}(r).
\end{equation}
Besides, as the expectation of $\mathbb{FE}_{t,s}$ is zero, the forecast is unbiased.
\end{proposition}
Theorem \ref{Invertibility1} entails that $\varepsilon_{r}\in\mathcal{M}_{s}$ for $r\le s$. Additionally, as $1\in\mathcal{M}_{s}$, it follows that $\sum_{r=-\infty }^{t}\xi (t,r)\varphi (r)\in\mathcal{M}_{s}$. Thus, both terms in the right-hand side of eq. (\ref{optimal linear predictor}) belong to $\mathcal{M}_{s}$. Using the formulas in eqs. (\ref{optimal linear predictor} and \ref{FE and MSE1}), one can conveniently verify that $\expectation\big(\mathbb{FE}_{t,s}\  \expectation(y_{t}|\mathcal{M}_{s})\big)=0$ (a detailed proof is provided in the Additional Appendix N). In other words $\mathbb{FE}_{t,s}$ and $\expectation(y_{t}|\mathcal{M}_{s})$ are orthogonal  (for short $\mathbb{FE}_{t,s} \perp \expectation(y_{t}|\mathcal{M}_{s})$), verifying that the predictor is optimal. 

Singh and Peiris (1987), Kowalski and Szynal (1990, 1991), and Grillenzoni (1990, 2000) obtained the evaluation of the optimal forecasts similar to that described in eq. (\ref{optimal linear predictor}), using prediction algorithms and a recursive procedure. We should also mention that Whittle (1965) showed that in general the linear least-square predictor obeys a recursion (see his eq. (12)) and gave a recursive method for obtaining its coefficients. In Section \ref{SECSTV} we describe how we can apply this technique in conjunction with our methodology to derive the optimal forecasts when the coefficients are varying stochastically. In particular, we examine the GRC-AR specification and a model with coefficients that follows AR processes.

An alternative and more realistic approach takes advantage of the explicit form of the process $y_t$, employed in the previous paragraph, but now given by eq. (\ref{TVAR(p)SOL1}). Let $\mathcal{K}_{s}$ be the smallest closed subspace of $L_2$ spanned by the set of $(p+q)$ past observations $\{y_{s},y_{s-1},\ldots y_{s+1-p}\}\cup \{\varepsilon_{s},\ldots,\varepsilon_{s+1-q}\}$, which also contains all constant functions. In a causal environment the $p$ past observations of 
the process $y_t$ in eq. (\ref{Wold Cramer}),  below and including the time point $s$, are realizations of the corresponding random variables, which are given by: $y_{s+1-m}=\sum_{r=-\infty}^{s+1-m}\xi (s+1-m,r)\upsilon_{r}$ for $1\leq m\leq p$.
Applying the conditional expectation operator based on  $\mathcal{K}_{s}$ to eq. (\ref{TVAR(p)SOL1}), the next  Proposition follows directly:

\begin{proposition}
\label{ProOptPred} The $k$-step-ahead optimal predictor of $y_t$  based on $\mathcal{K}_{s}$ is linear and is given by\vspace{-0.05in} 
\begin{equation*}
\expectation(y_{t}\left\vert \mathcal{K}_{s}\right. )=\sum_{m=1}^{p}\xi
^{(m)}(t,s)y_{s+1-m}+\sum_{r=s+1}^{t}\xi
(t,r)\varphi(r)+\sum\limits_{r=s+1-q}^{s}\xi _{s,q}(t,r)\varepsilon_{r}.
\end{equation*}%
The forecast error for the $k$-step-ahead predictor, that is $\mathbb{FE}_{t,s}=y_{t}-\expectation (y_{t}\left\vert\right. \mathcal{K}_{s})$ and its associated mean square error (its variance), coincide with those established by eqs. {\rm (\ref{FE and MSE1})}, that is:\vspace{-0.1in} 
\begin{equation}  \label{FE and MSE2}
\mathbb{FE}_{t,s}=\sum_{r=s+1}^{t}\xi _{q}(t,r)\varepsilon_{r},\ \mathbb{MSE}_{t,s}=\sum_{r=s+1}^{t}\xi _{q}^{2}(t,r)\mathbb{\sigma }^{2}(r).
\end{equation}
Besides, as the expectation of $\mathbb{FE}_{t,s}$ is zero, the forecast is unbiased.
\end{proposition}
In the Additional Appendix N, we show that $\expectation\big(\mathbb{FE}_{t,s}\ \expectation(y_{t}|\mathcal{K}_{s})\big)=0$  or equivalently that $\mathbb{FE}_{t,s} \perp \expectation(y_{t}|\mathcal{K}_{s})$, which is a prerequisite for the optimal nature of the predictors in $L_2$. As the forecast errors in Propositions \ref{Forecasting infinite observations} and \ref{ProOptPred} coincide (in $L_2$ sense), we conclude that:  $\expectation(y_{t}\left\vert \mathcal{K}_{s}\right.)=\expectation(y_{t}\left\vert \mathcal{M}_{s}\right.)$.
Moreover, in both cases the unconditional expectations and variances also coincide with that found earlier in Proposition \ref{ABSSUM}, that is:
$$\lim_{s\rightarrow -\infty }\expectation(y_{t}\left	\vert \mathcal{M}_{s}\right.)=\lim_{s\rightarrow -\infty }\expectation(y_{t}\left\vert \mathcal{K}_{s}\right.) =\sum \limits_{r=-\infty}^{t}\xi (t,r)\varphi(r)=\expectation (y_{t}),$$
$$\lim_{s\rightarrow -\infty }\mathbb{MS}\expectation_{t,s}=\sum_{r=-\infty }^{t}\xi _{q}^{2}(t,r)\sigma^{2}(r)=\var(y_{t}).$$ 
We remark that in a causal and invertible environment (the conditions of Theorem \ref{Invertibility1} hold) the past $q$ errors $\varepsilon_{r}$ for $s+1-q \le r\leq s$ can be expressed in terms of $y$'s by eq. (\ref{invertibility}). In practice,  eq. (\ref{invertibility}) can be used to provide error estimates for  $\{\varepsilon_{s},\ldots,\varepsilon_{s+1-q}\}$ by replacing in eq. (\ref{invertibility}) the $p$ realizations of the random variables: $\{y_{s},y_{s-1},\ldots y_{s+1-p}\}$. Furthermore, we remark that the explicit form of the variance for a Gaussian process is essential for the determination of the confidence intervals of $\expectation(y_{t}\left\vert \mathcal{K}_{s}\right.)$.

Finally, we formulate one of the main arguments made by Hallin (1986), which states that unlike the time invariant case, in a time-varying setting two
MSEs with the same forecasting horizon, but at different time points, are no longer equal. With this in mind, consider two pairs of time points, say $(t_{1},s_{1})$ and $(t_{2},s_{2})$, such that $t_{1}-s_{1}=t_{2}-s_{2}=k$.
The MSEs associated with these two time points are: 
\begin{equation*}
\mathbb{MSE}_{t_{1},s_{1}}=\sum_{r=s_{1}+1}^{t_{1}}\xi _{q}^{2}(t_{1},r)\mathbb{\sigma }^{2}(r),\ \ \ \mathbb{MSE}_{t_{2},s_{2}}=\sum_{r=s_{2}+1}^{t_{2}}\xi _{q}^{2}(t_{2},r)\mathbb{\sigma }^{2}(r).
\end{equation*}\\
Shifting the time intervals determined by the summation limits,  that is ($[s_{1}+1..t_{1}]$
and $[s_{2}+1..t_{2}]$), both to $[0..k-1]$, we get: 
\begin{equation*}
\mathbb{MSE}_{t_{1},s_{1}}=\sum_{r=0}^{k-1}\xi _{q}^{2}(t_{1},t_{1}-r)%
\mathbb{\sigma }^{2}(t_{1}-r),\ \ \ \mathbb{MSE}_{t_{2},s_{2}}=%
\sum_{r=0}^{k-1}\xi _{q}^{2}(t_{2},t_{2}-r)\mathbb{\sigma }^{2}(t_{2}-r).
\end{equation*}%
In a time-varying environment, a comparison between $\mathbb{MSE}%
_{t_{1},s_{1}}$ and $\mathbb{MSE}_{t_{2},s_{2}}$, whenever $%
t_{1}-s_{1}=t_{2}-s_{2}$, entails, in the general case, that : $\mathbb{MSE}%
_{t_1,s_1}\neq \mathbb{MSE}_{t_{2},s_{2}}$. On the other hand, in the time
invariant case $\xi _{q}(t,t-r)$, which will now be denoted as $\xi
_{r}^{(q)}$, becomes $\xi _{r}^{(q)}=\xi _{r}+\sum_{l=1}^{q}\xi _{r-l}\theta
_{l}$ (see eqs. (\ref{Toeplitz}) and (\ref{KSIq})). In this latter case the two MSEs coincide, as being identical to: 
\begin{equation*}
\mathbb{M}\mathbb{SE}_{k}=\sigma ^{2}\sum_{r=0}^{k-1}(\xi _{r}^{(q)})^{2}.
\end{equation*}

\subsection{Forward Asymptotic Efficiency\label{FORWASYMPTEFFIC}}

As pointed out by Hallin (1986), if a researcher wants to study the ``causal"
properties of the observed process, then he/she should examine the Wold-Cr%
\'{a}mer decomposition (see Theorem \ref{WoldCramDec} in Section \ref{WoldCramerDecomposition}). If forecasting is the main objective, then the
forecast obtained by the model should be asymptotically efficient in some
sense. Of course the asymptotic forecasting properties of the model rely on
its behaviour in the far future, whereas its causal properties involve its
remote past only. If the processes is of constant coefficients, these two issues coincide. In the case of time dependent coefficients, however, they apparently differ strongly.

To reiterate one of the main purposes in building models for stochastic processes is to provide convenient forecast procedures. The researcher would like to minimize (asymptotically) the MSE or in other words to achieve asymptotic efficiency. The asymptotic efficiency of a forecasting procedure can be defined in two alternative ways (seemingly, analogous to each other, but indeed basically different). The first one (termed by Hallin  efficiency, see Definition 5.1 in his paper) is obtained by considering the asymptotic forecasting performance of a model as $s$ tends to $-\infty$ (see Proposition \ref{ABSSUM}). The Wold-Cram\'{e}r representation  of a DTV-ARMA model in eq. (\ref{Wold Cramer}) is  efficient if and only if the model is invertible. If it is not an invertible model, then the mean square forecasting error associated with the above mentioned solution representation is not, in the general case, bounded. 

A more realistic approach to efficiency consists of considering the asymptotic behaviour of the mean square forecasting error as $t\rightarrow \infty $ for $s$ being arbitrary but fixed. This forward efficiency concept is also called the Granger-Andersen efficiency (see Definition 5.2 in Hallin, 1986, and the references therein). As the forecast error and the mean square error for the predictors given by eqs. (\ref{FE and MSE1} and \ref{FE and MSE2}), respectively, coincide, we can examine the forward asymptotic efficiency of DTV-ARMA models within the same forecasting scheme.

In the following Proposition, we give a weak condition, which guarantees the forward asymptotic boundedness of the mean square error.

\begin{proposition}
\label{prop.FAE} Let $F(t,s)\overset{\mathrm{def}}{=}\sum_{r=s+1}^{t}|%
\xi(t,r)|$ for $t>s$. If $\{F(t,s)\}_{t}$ is bounded, as a function of $t$
for each fixed $s$, then the mean square error is also bounded, as a
function of $t$. Equivalently, the boundedness of $\{F(t,s)\}_{t}$ entails
that for each $s$ either $\displaystyle \lim_{t\to \infty}\mathbb{MSE}_{t,s}$ exists in $\reals_{\ge 0}$ or $\{ \mathbb{MSE}_{t,s}\}_{t}$ oscillates with oscillation: $\mathit{\ \Omega}(s)=\displaystyle \inf_{t}(\sup_{r\ge t}\mathbb{MSE}_{r,s}-\inf_{r\ge t}\mathbb{MSE}_{r,s})$.
\end{proposition}
In the following Corollary, we provide a stronger condition, which guarantees the forward asymptotic uniform boundedness of the mean square error.
\begin{corollary}
\label{cor.FAE} If further the condition {\rm (\ref{abs summability condition2})} holds, that is the function $F(t)=\sum_{r=-\infty }^{t}|\xi(t,r)|$ is bounded in $t$, then $\mathbb{MSE}_{t,s}$ is uniformly bounded.
\end{corollary}
The proofs of the above statements are available in the Appendix.

\section{Stochastic Coefficients}\label{SECSTV}
In Sections \ref{Representations in L_2} and \ref{SecSecMom}, we restricted ourselves to the treatment of DTV models. In this Section we  examine processes with stochastically time
varying coefficients.  For simplicity instead of  ARMA processes, we will concentrate on the AR($p$) specification (STV-AR). In particular we will
investigate three models: the random coefficients one, its generalization and the double stochastic AR process. But first we will express the STV-AR
model in a companion matrix form.
The proofs of the present Section are provided in the Appendix.

\subsection{Companion Matrix Form}
In this Subsection we  show how to utilize the principal determinant (Green's function) and the $m$-th fundamental solution in order to
obtain a compact and explicit representation of the companion matrix.

The STV-AR($p$) process, can be expressed as 
\begin{equation}
y_{t}=\phi _{0t}+\boldsymbol{\phi }_{t}^{\prime }\mathbf{y}_{t-1}+\varepsilon_{t},  
  \label{STVAR}
\end{equation}%
where $\mathbf{y}_{t-1}=(y_{t-1,}$ $y_{t-2,}\ldots ,y_{t-p})^{\prime }$ is a 
$p\times 1$ vector of preceding random variables of $y_{t}$, and $\boldsymbol{\phi}_{t}=(\phi _{1t},$ $\phi _{2t},\ \ldots, \phi_{pt})^{\prime }$ is a $p\times 1$ vector of the autoregressive random coefficients. Notice that we
denote the STV coefficients, including the drift $\phi_{0t}$, by $\phi_{mt}$, $m=0,\ldots ,p$, instead of $\phi_{m}(t)$, which was the notation used for
the deterministic ones.

It is well known that the model (\ref{STVAR}) can be written in a companion
matrix form: 
\begin{equation}\label{COMPFORM}
\mathbf{y}_{t}=\boldsymbol{\phi}_{0t}+\mathbf{\Phi }_{t}\mathbf{y}_{t-1}+ \boldsymbol{\varepsilon}_{t},  
\end{equation}
where $\boldsymbol{\phi }_{0t}=(\phi_{0t}$ $0$ $\ldots 0)^{\prime }$, $\boldsymbol{\varepsilon }_{t}=(\varepsilon _{t}$ $0\ldots 0)^{\prime }$, and
the companion (square) matrix $\mathbf{\Phi }_{t}$ of order $p$ associated to the vector $\boldsymbol{\phi}_{t}$ is given by
\begin{equation} \label{COMPMATR}
\mathbf{\Phi }_{t}=\left( 
\begin{array}{ccccc}
\phi _{1t} & \phi _{2t} & \ldots & \phi _{p-1,t} & \phi _{pt} \\ 
1 & 0 & \ldots & 0 & 0 \\ 
0 & 1 & \ldots & 0 & 0 \\ 
\vdots & \vdots & \vdots \vdots \vdots & \vdots & \vdots \\ 
0 & 0 & \ldots & 1 & 0%
\end{array}%
\right) . 
\end{equation}%
That is, the STV-AR($p$) process is converted to a $p$-dimensional vector
STV-AR($1$) model. For any set of $p$ prescribed random variables $\mathbf{y}_{s}$, iterating eq. (\ref{COMPFORM}) yields 
\begin{equation}
\mathbf{y}_{t}=\mathbf{C}_{t,s}\mathbf{y}_{s}+\sum\limits_{r=s+1}^{t}\mathbf{C}_{t,r}(\boldsymbol{\phi}_{0,r}+\boldsymbol{\varepsilon }_{r}), 
\label{GSCOMPFORM}
\end{equation}%
where  $\mathbf{C}_{t,s}=\prod\limits_{r=s+1}^{t}\mathbf{\Phi}_{t,r}$ is the product of companion
matrices with initial square  matrix $\mathbf{C}_{t,t}=\mathbf{I}$ of order $p$. It follows directly from the above equation and Theorem \ref{TheoGenSol} (see also, for more details, Paraskevopoulos and Karanasos, 2021) that the $p$-dimensional square matrix $\mathbf{C}_{t,s}$ is given by 
\begin{equation*}
\mathbf{C}_{t,s}=\left(\begin{array}{cccc}
\xi^{(1)}(t,s) & \xi ^{(2)}(t,s) & \cdots & \xi^{(p)} (t,s) \\ 
\xi ^{(1)}(t-1,s-1) & \xi^{(2)}(t-1,s) & \cdots & \xi^{(p)}(t-1,s)
\\ 
\multicolumn{1}{c}{\vdots} & \multicolumn{1}{c}{\vdots} & \multicolumn{1}{c}{
\vdots \vdots \vdots} & \multicolumn{1}{c}{\vdots} \\ 
\xi^{(1)} (t-p+1,s) & \xi^{(2)} (t-p+1,s) & \cdots & \xi^{(p)}(t-p+1,s)
\end{array}\right).
\end{equation*}%
In other words the element occupying the $(n+1,m)$-th entry of the matrix 
$\mathbf{C}_{t,s}$ ($n=0,\ldots ,p-1$) is the $m$-th fundamental solution $\xi^{(m)}(t-n,s)$. We recall that $\xi^{(1)}(t,s)$ is given in eqs. (\ref{TVAR(P)1}
and \ref{ksi}) (where now, $\phi _{m}(t)$\ in eq. (\ref{TVAR(P)1}) is replaced by $\phi _{mt}$)
and, similarly, $\xi^{(m)}(t,s)$\ is given in eq. (\ref{Exdef:xi^m1}).

\subsection{Random Coefficients AR Model}
\label{Random Coefficients AR Model}
In this Subsection we  examine the random coefficient AR($p$) model (with acronym RC-AR($p$)), which is given by eq. (\ref{STVAR}), using for this the following notation and specifications: $\boldsymbol{\phi }_{t}^{\ast }=(\phi_{0t}$ $\boldsymbol{\phi }_{t})^{\prime }$,  $t=s+1,s+2\ldots $, is an $i.i.d.$ $(p+1)$-dimensional random vector of the coefficients, and the $i.i.d.$ errors, $\{\varepsilon_{t}, t\ge s+1\}$ are independent of the random drift and autoregressive coefficients. Let also $\varepsilon _{s}$ be a random variable independent of everything else, being the initial state.

Let us call $\left \vert \cdot \right \vert $ the Euclidean norm on the space $\reals^{p}$. Let $\reals^{p\times p}$ be the space of $p\times p$ matrices with elements in $\reals$ and denote by $\left \vert \left \vert \cdot \right \vert \right \vert $ the matrix norm induced by $\left \vert \cdot \right \vert $ (this is known as the spectral norm, defined as the largest singular value of the matrix).

\begin{condition}
\label{CONDRC}$\xi(t,s)\overset{a.s}{\rightarrow }0$ as $t\rightarrow \infty$.
\end{condition}
The proof of the next Theorem follows from the fact that $\xi _{t,s}\overset{a.s}{\rightarrow }0$ as $t\rightarrow \infty $ implies that $\left \vert \left \vert \mathbf{C}_{t,s}\right \vert \right \vert \overset{a.s.}{\rightarrow }0$ as $t\rightarrow \infty $ and, therefore, Theorem 2.1 in Erhardsson (2014) applies.
\begin{theorem}
\label{THEORRCAR}Consider the RC-AR($p$) model. Under Condition {\rm \ref{CONDRC}}
the following are equivalent:
\begin{center}
\begin{varwidth}{\textwidth}	
\begin{enumerate}
	\item[\it i)] $\mathbf{y}_{t}$ convergences in distribution as $t\to \infty$,
	\item[\it ii)]   $\sum \limits_{r=s+1}^{\infty }\left \vert \mathbf{C}_{r-1,s}\boldsymbol{\varepsilon}_{r}\right \vert <\infty$ a.s.,
	\item[\it iii)] $\sum \limits_{r=s+1}^{t}\mathbf{C}_{r-1,s}\boldsymbol{\varepsilon}_{r}$ converges a.s.,  as $t\to \infty$,
	\item[\it iv)] $\mathbf{C}_{t-1,s}\boldsymbol{\varepsilon }_{t}\overset{\rm a.s.}{\to}0$, as $t\to \infty$,
	\item[\it v)] $\sup_{t\ge s+1}\left\vert \mathbf{C}_{t-1,s}\boldsymbol{\varepsilon }_{t}\right\vert<\infty$ a.s.
\end{enumerate} 
\end{varwidth}
\end{center}
\end{theorem}
As pointed out by Erhardsson (2014) the implications $(ii) \implies (iii) \implies (iv) \implies (v)$ remain valid even if Condition (\ref{CONDRC}) does not hold.
\subsection{The Generalized RC-AR Model}
A generalization of the RC-AR process (GRC-AR($p$)) in eq. (\ref{STVAR}) is discussed in this Subsection. For this model, we show how the expected value of the $m$-th fundamental solution can be expressed as a Hessenbergian, which then can be employed to obtain an explicit representation of  the autocovariance function. 

The generalized process extends the simple one to cover contemporaneous dependence between the vector of the drift and the autoregressive coefficients, both included in $\boldsymbol{\phi}_{t}^{\ast}$, and the vector of errors $\varepsilon_{t}$. 

The GRC-AR($p$) model also integrates the following AR processes (for details, see Hwang and Basawa, 1998):
\begin{enumerate}
	\item[I.] Random coefficient model: $\phi _{mt}=\phi _{m}+\eta _{mt}$, where $\{\eta _{mt}\}$ is a sequence of $i.i.d.$ random variables with $\expectation(\eta_{mt})=0$, and $\{ \eta _{mt}\}$ is independent of $\{ \varepsilon_{\tau}, \tau\in\integers\}$ for all $m$ and $t\neq \tau$. Notice also that if we set $\eta_{mt}=0$, for all $m$ and $t$, we get the ordinary AR($p$) process.
	\item[II.] Markovian bilinear model: $\phi_{mt}=\phi_{m}+\vartheta_{m}\varepsilon_{t}$.
	\item[III.] Generalized Markovian bilinear model: $\phi _{mt}=c_{m}+\vartheta_{m}\varepsilon _{t}^{r_{m}}$, $r_{m}\in \integers_{>0}$, and $\varepsilon _{t}^{r_{m}}$ has finite moments $\expectation(\varepsilon_{t}^{r_{m}})$ for all $m$. Note that if we set $r_{m}=1$, for all $m$, we 	get the Markovian bilinear model coefficients.
    \item[IV.]	Random coefficient exponential model: $\phi_{mt}=c_{m}+(\vartheta_{m1}+\vartheta_{m2}e^{-\vartheta_{m3}\varepsilon_{t}^{2}})\varepsilon_{t}$.
\end{enumerate}
In cases I - IV all the coefficients without time $t$ subscript are constants.
In the next Proposition and Theorem we will make use of the following notation:

\begin{notation}\label{NOTExpVal}
\mbox{}	\vspace{-0.1in}\\
\begin{equation*}
\begin{array}{rcl}
y & \overset{\rm def}{=} & \expectation(y_{t}),\ \phi_{m}\overset{\rm def}{=}\expectation(\phi _{mt}),\ \sigma _{m\varepsilon }\overset{\rm def}{=} \expectation(\phi _{mt}\varepsilon _{t})\ \text{for}\ m=0, \ldots, p, \\ 
\sigma _{mn}&\overset{\rm def}{=}&\autocov(\phi _{mt},\phi _{nt}) = \bar{\phi}_{mn}-\phi_{m}\phi _{n}\ \text{ for } m,n=0, \ldots, p.
\end{array}%
\end{equation*}
\end{notation}

\begin{condition}
\label{CONDGRCARM1}$\left \vert \sum \nolimits_{m=1}^{p}\phi
_{m}\right
\vert <1.$
\end{condition}

\begin{proposition}
\label{APPROPGRCMEAN}Under Condition {\rm \ref{CONDGRCARM1}} the expected value of
the GRC-AR($p$) model is given by 
\begin{equation*}
y=\frac{\phi_{0}}{1-\sum \limits_{m=1}^{p}\phi _{m}}.
\end{equation*}
\end{proposition}

In what follows we adopt the notation $\mathbf{\Phi}_{t}^{\otimes 2}=\boldsymbol{\Phi }_{t}\otimes\mathbf{\Phi }_{t}$, where $\boldsymbol{\Phi }_{t}$ is the companion matrix   (see eq. (\ref{COMPMATR})), $\otimes $ is the Kronecker product, and $\lambda _{\max }[\expectation(\mathbf{\Phi }_{t}^{\otimes 2})]$ refer to the modulus of the largest eigenvalue of $\expectation(\mathbf{\Phi }_{t}^{\otimes 2})$.

\begin{condition}
\label{CONDGRCARM2}$\lambda _{\max }[\expectation(\mathbf{\Phi }_{t}^{\otimes
2})]<1$.
\end{condition}

In the next Theorem we shall also make use of the following notation:

\begin{notation}\label{NOTGRC} $\mbox{}$\newline
\mbox{}\\
i) $\varphi _{ij}$ denotes the element occuping the $(i,j)$-th entry of the
matrix \ $[\mathbf{I}_{p^{2}}-\mathbb{E(}\mathbf{\Phi }_{t}^{\otimes
2})]^{-1}$, 
\begin{equation}\label{SIGMA}
\hspace{-2.5in} ii)\text{ }\sigma ^{2}=(\sigma _{\varepsilon }^{2}+\sigma
_{00}+2\sigma _{0\varepsilon })+2y\sum\limits_{m=1}^{p}(\sigma _{m\varepsilon
}+\sigma _{m0})+y^2\sum\limits_{n=1}^{p}\sum\limits_{m=1}^{p}\sigma _{mn},
\end{equation}
where $\sigma _{\varepsilon }^{2}=\var(\varepsilon _{t})$. 
\hspace{-0.3in} iii) $\xi _{k}^{(m)}=\expectation(\xi^{(m)}(t,s))$, recall that $k=t-s$, is the $k$-dimensional banded Hessenbergian: 
\begin{equation}
\xi _{k}^{(m)}=\left\vert 
\begin{array}{cccccccc}
{\footnotesize \phi }_{m} & {\footnotesize -1} &  &  &  &  &  &  \\ 
{\footnotesize \phi }_{m+1} & {\footnotesize \phi }_{1} & \ddots &  &  &  &  & 
\\ 
{\footnotesize \vdots } & {\footnotesize \vdots } & \ddots & \ddots &  &  & 
&  \\ 
{\footnotesize \phi }_{p} & {\footnotesize \phi }_{p-1} & \ddots & \ddots & 
\ddots &  &  &  \\ 
& {\footnotesize \phi }_{p} & \ddots & \ddots & \ddots & \ddots &  &  \\ 
&  & {\small \ddots } & {\small \ddots } & {\small \ddots } & {\small \ddots 
} & \ddots &  \\ 
&  &  & {\footnotesize \phi }_{p} & {\footnotesize \phi }_{p-1} & {\small %
\cdots } & {\footnotesize \phi }_{1} & {\footnotesize -1} \\ 
&  &  &  & {\footnotesize \phi }_{p} & {\small \cdots } & {\footnotesize %
\phi }_{2} & {\footnotesize \phi }_{1}%
\end{array}%
\right\vert  \label{KSIMEAN}
\end{equation}%
(we recall that $\xi^{(1)}(t,s)$\ is given in eqs. {\rm(\ref{TVAR(P)1}} and {\rm \ref{ksi})}, where now, $\phi_{m}(t)$\ is replaced by $\phi _{mt}$ and, similarly, $\xi^{(m)} (t,s)$ is given in eq.  {\rm (\ref{Exdef:xi^m1}))}.
\end{notation}

\begin{theorem}
\label{APPTHGRCCOV}\bigskip Under Condition {\rm \ref{CONDGRCARM2}} the covariance
structure of the GRC-AR($p$) model is given by 
\begin{equation*}
\gamma (\ell )=\left\{ 
\begin{array}{ccc}
\varphi_{\ell +1,1}\sigma^{2} & \text{for} & \ell =0,\ldots ,p-1 \\ 
\sum\limits_{m=1}^{p}\xi _{\ell }^{(m)}\gamma (m-1) & \text{for} & \ell \geq
p, \end{array}\right.
\end{equation*}
where $\xi_{\ell }^{(m)}$ is given by eq. {\rm (\ref{KSIMEAN})}.
\end{theorem}
Let us remark that Condition \ref{CONDGRCARM1} is equivalent to  $\lim_{\ell \rightarrow \infty }\xi _{\ell }^{(1)}=0$, which implies that $\lim_{\ell \rightarrow \infty }\xi _{\ell }^{(m)}=0$ for all $m$, and, therefore, ensures that $\lim_{\ell \rightarrow \infty }\gamma (\ell )=0$. As another example we examine the autoregressive process of order 2.

\begin{corollary}
Consider the GRC-AR($2$) model and let the following condition hold: 
\begin{equation*}
1-\frac{2\bar{\phi} _{12}\phi _{1}}{1-\phi _{2}}>\bar{\phi} _{11}+\bar{\phi} _{22}>0.
\end{equation*}%
Under the above condition its covariance structure is given by 
\begin{equation*}
\ \gamma (\ell )=\left \{ 
\begin{array}{ccc}
\varphi _{\ell +1,1}\sigma ^{2} & \text{for} & \ell =0,1 \\ 
\sum \limits_{m=1}^{2}\xi _{\ell }^{(m)}\gamma (m-1) & \text{for} & \ell
\geq 2,%
\end{array}%
\right.
\end{equation*}%
where 
\begin{eqnarray*}
\varphi _{1,1} &=&\frac{1-\phi _{2}^{2}}{(1-\phi _{2}^{2})(1-\bar{\phi} _{11}-\bar{\phi}
_{22})-2\bar{\phi} _{12}\phi _{1}(1+\phi _{2})}, \\
\varphi _{2,1} &=&\varphi _{2,1}=\frac{\phi _{1}(1+\phi _{2})}{(1-\phi
_{2}^{2})(1-\bar{\phi} _{11}-\bar{\phi }_{22})-2\bar{\phi} _{12}\phi _{1}(1+\phi _{2})},
\end{eqnarray*}%
and the tridiagonal matrices of order $\ell $, that is $\xi_{\ell}^{(m)}$, 
$m=1,2$, are given by 
\begin{equation*}
\xi _{\ell }^{(1)}=\left[ 
\begin{array}{ccccc}
{\footnotesize \phi }_{1} & {\footnotesize -1} &  &  &  \\ 
{\footnotesize \phi }_{2} & {\footnotesize \phi }_{1} & \ddots &  &  \\ 
& \ddots & \ddots & \ddots &  \\ 
&  & \ddots & {\footnotesize \phi }_{1} & {\footnotesize -1} \\ 
&  &  & {\footnotesize \phi }_{2} & {\footnotesize \phi }_{1}%
\end{array}%
\right] ,\text{ \ }\xi _{\ell}^{(2)}={\footnotesize \phi }_{2}\xi_{\ell -1}^{(1)}.
\end{equation*}
\end{corollary}

\subsection{Double Stochastic AR Models}
\label{Double Stochastic AR Models}
In this Subsection we examine the more general case where the autoregressive coefficients follow AR processes. We  show that for this
model the unconditional variance exists in $\reals_{>0}$ provided that the associated Green function convergences in $L_{2}$, a result which is in line with Theorem 2(ii). In other words, we investigate the double stochastic AR model, hereafter termed DS-AR (for double stochastic processes, and in particular ARMA processes with ARMA coefficients see, for example, Grillenzoni, 1993, and the references therein). The DS-AR model is defined by eq. (\ref{STVAR}) but in this case the autoregressive coefficients, $\phi_{mt}$ for $m=1,\ldots ,p$, follow  AR processes: 
\begin{equation}\label{DARcoef}
\phi _{mt}=\beta_{m0}+\sum \limits_{l=1}^{p_{m}}\beta _{ml}\phi_{m,t-l}+e_{mt},  
\end{equation}
where $\beta_{m0}$ and $\beta_{ml}$ are constant coefficients and $p_{m}\in \integers_{> 0}$ for all $m:  1\le m\le p$. $\{e_{mt}\}$ are martingale difference sequences  defined on $L_{2}$, where $e_{mt}$ and $\varepsilon _{t\pm b}$, $b\in \integers$, are independent of each other for all $m$, and $t\in\integers$. For simplicity, we will assume that the drift in eq. (\ref{STVAR}) is time invariant, that is, $\phi _{0t}=\phi _{0}$ for all $t$.

The results in Sections \ref{Representations in L_2} and \ref{SecSecMom} can be easily modified to cover DS-AR models by replacing the fundamental solutions with their respective (unconditional and conditional) expectations. More specifically, we present two Theorems followed by two Propositions (their proofs essentially repeat the arguments of the proofs of those in Section \ref{Representations in L_2} and \ref{SecSecMom}).

Sufficient conditions ensuring the Wold-Cram\'{e}r decomposition of DS-AR($p$) models, and therefore the existence of the first two unconditional moments, are the following two:
\begin{equation}
\left. 
\begin{array}{lll}
\displaystyle\sum_{r=-\infty }^{t}|\expectation\big(\xi (t,r)\big)| & < & \infty ,\ 
\text{for\ all}\ t\ \text{(first-order\ absolute\ summability)}, \\ 
\displaystyle\sum_{r=-\infty }^{t}\expectation\big(\xi ^{2}(t,r)\big) & < & \infty ,\ 
\text{for\ all}\ t\ \text{ (second-order\ summability).}%
\end{array}
\right\}  \label{2ndORDSUMMABILITY}
\end{equation}

\begin{remark}
Generally, it is very difficult to verify if the two summability conditions are fulfilled. Only some special cases allow to write explicit solutions
(see,  And\u{e}l, 1991, and the references therein). A sufficient condition for the absolute summability to hold is that $\{ \sum \nolimits_{m=1}^{p}\phi
_{mt}\}$ belongs with probability one to the interval $(-1,1)$, nearly  everywhere, that is, with the exception, at most, of a finite number of $t$
(see, for example, Grillenzoni, 1993). Similarly, a sufficient condition for the square summability to hold is that with probability one  $\lambda _{t}^{(\max )}[\mathbf{\Phi}_{t}^{\otimes 2}]<1$, nearly everywhere, where $\lambda _{t}^{(\max )}[\mathbf{\Phi }_{t}^{\otimes 2}]$ refers to the modulus of the largest eigenvalue of $\mathbf{\Phi }_{t}^{\otimes 2}$; we recall that $\mathbf{\Phi }_{t}$ is the companion matrix (see eq. (\ref{COMPMATR})).
\end{remark}

\begin{theorem}
\label{StoProWoldRepr} Let the two summability conditions in {\rm(\ref{2ndORDSUMMABILITY})} hold. The Wold-Cram\'{e}r decomposition (in $L_2$
sense) is a solution of the DS-AR($p$) model in  eq.  {\rm (\ref{STVAR})}, where its autoregressive  coefficients are given by eq. {\rm(\ref{DARcoef})}, being of the form: 
\begin{equation}\label{WoldCramer1}
y_{t}=\sum_{r=-\infty }^{t}\xi (t,r)(\phi _{0}+\varepsilon _{r}).
\end{equation}
\end{theorem}

\begin{theorem}
\label{THEORDSARWC}If all the autoregressive coefficients, $\phi _{mt}$, $m=1,\ldots ,p$, are strictly stationary, then eq. {\rm(\ref{STVAR})} has a
stationary solution of the type {\rm(\ref{WoldCramer1})} if and only if 
\begin{equation*}
\sum_{r=1}^{\infty }\left \vert \expectation\big(\xi (r,1)\big)\right \vert <\infty \  \ 
\text{\rm and} \ \  \  \sum_{r=1}^{\infty }\expectation\big(\xi ^{2}(r,1)\big)<\infty .
\end{equation*}
\end{theorem}

In what follows, we present explicit formulae for the first and second unconditional moments for the DS-AR family of processes coupled with sufficient and necessary conditions for their existence.

\begin{proposition}
\label{PROPABSSUM1} Under the two conditions in {\rm(\ref{2ndORDSUMMABILITY})},
it follows from Theorem {\rm\ref{StoProWoldRepr}} that the unconditional mean of the DS-AR($p$) process, $y_{t}$, exists in $\reals$ and is given by 
\begin{equation}
\expectation(y_{t})=\phi _{0}\sum\limits_{r=-\infty }^{t}\expectation\big(\xi (t,r)\big).
\label{first unconditional moment1}
\end{equation}%
A necessary condition for the absolute summability to hold is 
\begin{equation*}
\lim_{s\rightarrow -\infty }\expectation\big(\xi (t,s)\big)=0.
\end{equation*}
Moreover, the unconditional variance of the process exists in $\reals_{> 0}$ and it is given by 
\begin{equation*}
\var(y_{t})=\phi _{0}^{2}\sum\limits_{r=-\infty }^{t}\mathbb{V}ar\big(\xi (t,r)\big)+\sigma_{\varepsilon }^{2}\sum_{r=-\infty}^{t}\expectation\big(\xi
^{2}(t,r)\big).
\end{equation*}%
A necessary condition for the second-order summability to hold is 
\begin{equation*}
\lim_{s\rightarrow -\infty }\expectation\big(\xi ^{2}(t,s)\big)=0.
\end{equation*}
\end{proposition}

We notice that $\lim_{s\rightarrow -\infty }\expectation\big(\xi^{2}(t,s)\big)=0$ is
equivalent to $\lim_{s\to -\infty}||\xi (t,s)||^2_{L_2}=0$, which, in turn,
is equivalent to $\lim_{s\to -\infty}||\xi (t,s)||_{L_2}=0$. In this latter
case, we write $\lim_{s\to -\infty}\xi (t,s)\overset{L_2}{=}0$, which is in line with Theorem \ref{theo: asympt. stability in L2}(ii).

\begin{proposition}\label{PROBOPTPRED}
Following the notation of Proposition {\rm\ref{ProOptPred}}, let $\mathcal{K}_{s}$ be the smallest closed subspace of $L_2$ spanned by the finite observed sequence\vspace{-0.1in} 
$$\{y_s,...,y_{s+1-p}\}\cup \bigg(\displaystyle\bigcup_{i=1}^m\{\phi _{i,s},\phi_{i,s+1},\ \ldots, \phi_{i,s+1-p_i}\}\bigg).$$
The $k$-step-ahead optimal (in $L_{2}$-sense)
linear predictor of the DS-AR($p$) process is:
\begin{equation*}
\expectation(y_{t}\left\vert \mathcal{K}_{s}\right. )=\sum_{m=1}^{p}\expectation\big(\xi ^{(m)}(t,s)|\mathcal{K}_{s}\big)y_{s+1-m}+\phi _{0}\sum_{r=s+1}^{t}\expectation\big(\xi (t,r)\left\vert \mathcal{K}_{s}\right.\big).
\end{equation*}%
We remark that $\lim_{s\to -\infty}\expectation(y_{t}\left\vert \mathcal{K}_{s}\right.)=\expectation(y_{t})$ is  given by eq. (\ref{first unconditional moment1}).
In addition, the forecast error for the above $k$-step-ahead predictor, $\mathbb{FE}_{t,s}$, is given by 
\begin{equation*}
\mathbb{FE}_{t,s}=\phi _{0}\sum_{r=s+1}^{t}\bigg(\xi (t,r)-\expectation(\xi
(t,r)\left\vert \mathcal{K}_{s}\right.)\bigg)+  \sum_{r=s+1}^{t}\xi(t,r)\varepsilon_{r}+ 
\sum_{m=1}^{p}\big(\xi ^{(m)}(t,s)-\expectation(\xi
^{(m)}(t,s)|\mathcal{K}_{s})\big)y_{s+1-m}.\ 
\end{equation*}%
The  conditional variance of $y_t$ based on $\mathcal{K}_{s}$, that is $\var(y_t
\left\vert \mathcal{K}_{s}\right. )=\var(\mathbb{FE}
_{t,s}\left\vert \mathcal{K}_{s}\right. )$ is given by 
\begin{equation*}
\var(y_t
\left\vert \mathcal{K}_{s}\right. )=\phi _{0}^{2}\sum_{r=s+1}^{t}\var\big(\xi (t,r)\left\vert 
\mathcal{K}_{s}\right.\big)+\sigma_{\varepsilon }^{2}\sum_{r=s+1}^{t} \expectation\big(\xi ^{2}(t,r)\left\vert \mathcal{K}_{s}\right.\big)+\sum_{m=1}^{p}\var\big(\xi^{(m)}(t,s)\left\vert \mathcal{K}_{s}\right. \big)y_{s+1-m}^{2}
\end{equation*}%
(we recall that $\xi (t,r)$ have been introduced in eqs. {\rm(\ref{TVAR(P)1})} and {\rm (\ref{ksi})}, where now $\phi _{m}(t)$ are replaced by $\phi_{mt}$). \\
We remark that $\lim_{s\to -\infty}\var (y_{t}\left\vert \mathcal{K}_{s}\right.)=\var(y_{t})$, which is given in Proposition \ref{PROPABSSUM1}.
\end{proposition}

\section{Time-varying Polynomials\label{SECTVOPERATORS}}

In Sections \ref{SubsecMainTheor}, \ref{Representations in L_2} and \ref{SecSecMom},  we employed techniques of linear algebra in order to obtain an explicit  representation of the TV-ARMA($p,q$) model and its
first two (conditional and unconditional) moments. The main mathematical tool used was the Hessenbergian determinant. Now that we have expressed the
Green's function as a Hessenbergian, we will see how the summation terms in the various equations in Sections \ref{SubsecMainTheor} and \ref{Representations in L_2} can be expressed as time-varying polynomials.

Recall that $B$\ denotes the backshift (or lag operator), defined such that $%
By_{t}=y_{t-1}$. The time-varying AR and moving average (MA) polynomial (backshift) operators associated with the TV-ARMA($p,q$) model are denoted
as:\vspace{-0.1in}
\begin{equation}
\mathbf{\Phi}_{t}(B)=1-\sum \limits_{m=1}^{p}\phi_{m}(t)B^{m},\text{ }\Theta
_{t}(B)=1+\sum \limits_{l=1}^{q}\theta_{l}(t)B^{l}.  \label{operators}
\end{equation}
Under this notation eq. (\ref{TVAR(P)}) can be written in a more condensed
form%
\begin{equation}
\Phi_{t}(B)y_{t}=\varphi(t)+\Theta_{t}(B)\varepsilon_{t}.  \label{TVAR(p)2}
\end{equation}

In the time invariant case one can employ the roots of the time invariant polynomial $\Phi (z^{-1})$\ to obtain its general time series properties such as the Wold decomposition and the second moment structure. In a time-varying environment, according to Grillenzoni (1990), the so called here principal fundamental sequence $\{\xi (t,s)\}_{t\geq s+1-p}$ cannot be obtained as in stationarity, that is by expanding in Taylor series the rational polynomial $\Phi _{t}^{-1}(B)$. As an alternative, Hallin (1986) used some results on difference operators involving the symbolic product of operators, which has also been termed by researchers, in the field of engineering, as the skew multiplication operator (see, for example, Mrad and Farag, 2002). Hence, now that we have at our disposal an explicit and computationally tractable representation of the Green function as a banded Hessenbergian, coupled with the use of the time-domain noncommutative multiplication operation -which, as pointed out by Mrad and Farag (2002), is based on the manipulation of polynomial operators with time-varying coefficients using operations restricted to the time domain- we are able to state some important Theorems in relation to the results in Sections \ref{SubsecMainTheor} and \ref{Representations in L_2}.

\subsection{The Skew Multiplication Operator\label{SUBSKEW}}

In a time-varying environment, the time-varying polynomial operators in eq. (\ref{operators}) can be manipulated by using the skew multiplication
operator \textquotedblleft$\circ$\textquotedblright \ defined by
\begin{equation}
B^{i}\circ B^{j}=B^{i+j}\text{ \ and \ }B^{i}\circ f(t)=f(t-i)B^{i}\text{,}
\label{skew}
\end{equation}
where $f(t)$ is a function of time. This time-domain multiplication
operation is associative but noncummutative (see Karmen, 1988, Bouthellier
and Ghosh, 1988, and Mrad and Farag, 2002). Using the properties of
\textquotedblleft$\circ$\textquotedblright, from eq. (\ref{TVAR(p)2}), under
the necessary and sufficient conditions in Proposition \ref{ABSSUM}, we can
obtain the unique inverse of $\Phi_{t}(B)$, that is $\Phi^{-1}_{t}(B)\circ
\Phi_{t}(B)=1$, where $\Phi^{-1}_{t}(B)$ is provided in eq. (\ref{KSIUncond}) below.

\subsection{Polynomial Operators}

Next, and equally important, we will provide a critical and essential further link between the linear algebra techniques used in Section \ref{SubsecMainTheor} (to obtain the explicit representation of the TV-ARMA model in eq. (\ref{TVAR(p)SOL1})) and the time-varying polynomial approach, in which we make use of the skew multiplication operator. Certainly, from an operational point of view, both are equally satisfying and recommendable. First, let us define the two time-varying polynomial (backshift) operators associated with eq. (\ref{TVAR(p)SOL1}).

\begin{definition}
\label{DEFKSIHOMPAR} I) Let $\Xi _{t,p}^{[k]}(B)$ be defined as follows:
\begin{equation}\label{KSIHOM}
\Xi _{t,p}^{[k]}(B)=1-\sum \limits_{m=1}^{p}\xi ^{(m)}(t,s)B^{k-1+m}.
\end{equation}

IIa) Let $\Xi _{t}^{(k)}(B)$ be defined as follows:
\begin{equation}\label{KSI}
\Xi _{t}^{(k)}(B)=\sum_{r=s+1}^{t}\xi (t,r)B^{t-r}\overset{\mathrm{or}}{=}
\sum_{r=0}^{k-1}\xi (t,t-r)B^{r}. 
\end{equation}

IIb) The limit of $\Xi _{t}^{(k)}(B)$ as $k\rightarrow \infty $ is denoted
by $\Xi _{t}(B)$.
\end{definition}
Notice that we have employed the notation $(.)^{[k]}$ ($k$ inside brackets) for the polynomial operator defined in eq. (\ref{KSIHOM}) in order to be distinguished from the one defined in eq. (\ref{KSI}).
\begin{remark}
\label{REFKSITPK}$\Xi _{t,p}^{[k]}(B)$ in eq. {\rm (\ref{KSIHOM})} is a polynomial of order $k-1+p$ associated with the homogeneous solution in eq.  {\rm(\ref{alternative homogeneous solution})}, and it is expressed in terms of the 
$m$ fundamental solutions, defined in eq. {\rm(\ref{Exdef:xi^m1})}.\newline
Notice that: i) $\Xi _{t,p}^{[1]}(B)=\Phi _{t}(B)$, since $\xi^{(m)}(t,t-1)=\phi _{m}(t)$ {\rm (}see the discussion next to eq. {\rm (\ref{Exdef:xi^m1}))}, and\\ ii) Under the stability condition in Theorem \ref{theo: asympt. stability in L2}(i), that is for the deterministic case, $\lim_{k\rightarrow \infty }\Xi _{t,p}^{[k]}(B)=1$, since $\lim_{s\rightarrow -\infty }\xi ^{(m)}(t,s)=0$.
\end{remark}

\begin{remark}
$\Xi _{t}^{(k)}(B)$ is a polynomial of order $k-1$ associated with the particular solution in eq. {\rm(\ref{PARTICULAR SOLUION})}, and it is expressed in terms of the first $k$ instances of the principal determinant, $\xi(t,s+r)$ for $1\le r\le k$. Notice also that $\Xi _{t}^{(1)}(B)=1$.
\end{remark}

Next we define two additional time-varying polynomial operators associated
with the innovation part of the particular solution (see eq. \ref{eq. condensed}) and the Wold-Cram\'{e}r decomposition (see eq. \ref{Wold Cramer}) respectively.

\begin{definition}
\label{KSIq,t,s} I) Let $\Xi _{t,q}^{(k)}(B)$ be defined as follows%
\begin{equation}\label{KSIQTS}
\Xi _{t,q}^{(k)}(B)=\sum_{r=0}^{k-1}\xi
_{q}(t,t-r)B^{r}+\sum_{r=k}^{k-1+q}\xi _{s,q}(t,t-r)B^{r}.  
\end{equation}%
II) Let $\Xi _{t,q}(B)$ be defined as follows 
\begin{equation*}
\Xi _{t,q}(B)=\sum_{r=0}^{\infty }\xi _{q}(t,t-r)B^{r},
\end{equation*}
that is:  $\lim_{k\to \infty}\Xi _{t,q}^{(k)}(B)=\Xi _{t,q}(B)$.
\end{definition}

\begin{remark}
\label{REMKSITQK} $\Xi _{t,q}^{(k)}(B)$ is a polynomial of order $k-1+q$
associated with the innovation part of the particular solution, and is expressed in terms of the first $k$ instances of  $\xi_{q}(t,s+r)$ for $1\le r\le k$ and the first $q$ instances of $\xi_{s,q}(t,s+r-q)$ for $1\le r\le q$, which have been introduced in Definition {\rm\ref{DEFKSIQ}}.\newline
Notice that, i) for the pure AR model $\Xi_{t,0}^{(k)}(B)=\Xi _{t}^{(k)}(B)$, since $\xi _{0}(t,t-r)=\xi (t,t-r)$ and the second summation in eq. {\rm(\ref{KSIQTS})} (adopting the convention $\sum_{r=k}^{k-1}a_r=0$) vanishes. ii) $\Xi _{t,q}^{(1)}(B)=\Theta _{t}(B)$, since the first summation is equal to $\xi (t,t)=1$, and the second summation is equal to $1-\Theta _{t}(B)$ (see
the discussion next to Proposition \ref{PROPDECOMP}).
\end{remark}

\subsection{Polynomial Representations}
The next Proposition is an alternative to Proposition \ref{PROPDECOMP}, expressed in terms of polynomial operators. Proposition \ref{PROPGen2b} and Theorems \ref{THEORGENERPOLYNOM}, \ref{TheoremWoldPolynomials} below can be deduced by applying the properties of the skew multiplication operator \textquotedblleft $\circ $\textquotedblright (see eq. (\ref{skew})). By analogy to the decomposition in eq. (\ref{eq. condensed}), we deduce the following two part-decomposition:

\begin{proposition}
\label{PROPGen2b} $\Xi_{t}^{(k)}(B)\circ u_{t}$ takes the following
alternative expressions: 
\begin{equation}
\Xi_{t}^{(k)}(B)\circ u_{t}= \Xi _{t,q}^{(k)}(B)\varepsilon_{t}\overset{%
\mathrm{or}}{=}\sum_{r=0}^{k-1}\xi_{q}(t,t-r)\varepsilon_{t-r}+%
\sum_{r=k}^{k-1+q}\xi _{s,q}(t,t-r)\varepsilon_{t-r}.  \label{Decompose}
\end{equation}
\end{proposition}
A proof of the above result is provided in Appendix.
The first of the following Theorems (Theorem \ref{THEORGENERPOLYNOM}) is
equivalent to Theorem \ref{TheoGenSol}.

\begin{theorem}
\label{THEORGENERPOLYNOM} The explicit  representation of eq. {\rm(\ref{TVAR(P)})} given
by eq. {\rm(\ref{TVAR(p)SOL1})}, can be equivalently expressed, in terms of the
polynomial operators given in Definitions {\rm\ref{DEFKSIHOMPAR}} and {\rm\ref{KSIq,t,s}}, as: 
\begin{equation}
\Xi _{t,p}^{[k]}(B)y_{t}=\Xi _{t}^{(k)}(B)\circ \lbrack \varphi
(t)+u_{t}]=\Xi _{t}^{(k)}(B)\circ \varphi (t)+\Xi _{t,q}^{(k)}(B)\varepsilon
_{t}.  \label{Gen2a}
\end{equation}
\end{theorem}
Clearly, the results for the pure AR model are obtained by setting $u_{t}=\varepsilon _{t}$ in the first equality of eq. (\ref{Gen2a}) or by
noticing that $\Xi _{t}^{(k)}(B)=\Xi _{t,0}^{(k)}(B)$ in the second equality of eq. (\ref{Gen2a}) (see Definition \ref{KSIq,t,s} and Remark \ref{REMKSITQK}). Notice also that
Proposition \ref{ProOptPred} can be expressed in terms of the time-varying
polynomial operators as well (results not reported here).

In what follows we will make use of the infinite order polynomials, $\Xi
_{t}(B)$ and $\Xi _{t,q}(B)$, which have been introduced in Definitions \ref%
{DEFKSIHOMPAR}(IIb) and \ref{KSIq,t,s}(II), respectively. The next Theorem shows that $\Xi _{t,q}(B)$ is the time-varying Wold-Cram\'{e}r polynomial operator associated with the DTV-ARMA($p,q$) model.

\begin{theorem}
\label{TheoremWoldPolynomials} Under the absolute summability condition in {\rm (\ref{abs summability condition})}, the Wold-Cram\'{e}r operators are given by 
\begin{equation}\label{KSIUncond}
\Xi _{t}(B)=\Phi_{t}^{-1}(B)\ \  \text{\rm and }\ \   \Xi _{t,q}(B)=\Phi _{t}^{-1}(B)\circ \Theta _{t}(B)\ \ \
\  \mathrm{w.r.t}\ \ ``\circ ".
\end{equation}
The Wold-Cram\'{e}r decomposition in Theorem {\rm\ref{WoldCramDec}} can be written in terms of the Wold-Cram\'{e}r operator as 
\begin{equation}\label{WoldAgain}
y_{t}=\Xi_{t}(B)\circ \varphi(t)+\Xi_{t,q}(B)\varepsilon_{t},
\end{equation}%
which of course implies that
\begin{equation*}
\mathbb{E(}y_{t})=\Xi_{t}(B)\circ \varphi(t).
\end{equation*}
\end{theorem}
Appendix Section \ref{APPENDTVARPOLYNOM} contains the proofs of both Theorems \ref{THEORGENERPOLYNOM} and \ref{TheoremWoldPolynomials}. In the
Online Appendix I, we show how the time-varying techniques, introduced here, that is, polynomials with the time dependent coefficients (time-varying polynomials for short) expressed as Hessenbergians (coupled with the usage of the skew multiplication operator), incorporate as a special case the standard approach to time series analysis which is based on characteristic polynomials. In the Online Appendix J, we present a summary of the time-varying polynomials results.

\section{An Example}\label{Example}

In this Section we will consider the AR($2$) process with $2$ deterministic abrupt breaks [DAB-AR($2;2$)] at fixed points of time $t_{1}$ and $t_{2}$,
where $t_{1}>t_{2}$.

\begin{example}
\label{ex. deterministic abrupt breaks} The process is defined by 
\begin{equation}
y_{t}=\left \{ 
\begin{array}{ccc}
\varphi _{1}+\phi _{1,1}y_{t-1}+\phi _{2,1}y_{t-2}+\sigma _{1}e_{t} & \text{%
for} & t>t_{1}, \\ 
\varphi _{2}+\phi _{1,2}y_{t-1}+\phi _{2,2}y_{t-2}+\sigma _{2}e_{t} & \text{%
for} & t_{2}<t\leq t_{1}, \\ 
\varphi _{3}+\phi _{1,3}y_{t-1}+\phi _{2,3}y_{t-2}+\sigma _{3}e_{t} & \text{for} & t\leq t_{2},%
\end{array}%
\right.  \label{AB-AR(2,2)}
\end{equation}%
where $e_{t}\sim $ i.i.d. $(0,1)$ for all $t$ and $0<\sigma _{i}^{2}\leq M$, 
$i=1,2,3$. Applying the results of Theorem {\rm\ref{TheoGenSol}} the
following Corollary provides the solution representation  (at time $t_{1}+l$, $l\in \integers_{\geq 0}$) of eq. {\rm (\ref{AB-AR(2,2)})}. But first, we define the following matrices.
\end{example}

\begin{definition}
Let the two tridiagonal matrices of order $r\in \integers_{\geq 1}$,
denoted by $\Phi _{r}^{(j)}$, $j=1,2$, be defined as%
\begin{equation*}
\mathbf{\Phi }_{r}^{(j)}=\left[ 
\begin{array}{cccccc}
\phi _{1,j} & -1 &  &  &  &  \\ 
\phi _{2,j} & \phi _{1,j} & -1 &  &  &  \\ 
& \phi _{2,j} & \phi _{1,j} & -1 &  &  \\ 
&  & \ddots & \ddots & \ddots &  \\ 
&  &  & \phi _{2,j} & \phi _{1,j} & -1 \\ 
&  &  &  & \phi _{2,j} & \phi _{1,j}%
\end{array}%
\right] ,
\end{equation*}
\end{definition}

\begin{definition}
i) The tridiagonal matrix (of order $l-r$) $\mathbf{\Phi} _{t_{1}+l,t_{1}+r}$, for $%
r=1,\ldots ,l-1$ and $l\geq 1$, is defined as: 
\begin{equation*}
\mathbf{\Phi} _{t_{1}+l,t_{1}+r}=\mathbf{\Phi} _{l-r}^{(1)},
\end{equation*}
where its determinant is $\xi (t_{1}+l,t_{1}+r)=\left \vert \mathbf{\Phi}
_{t_{1}+l,t_{1}+r}\right \vert $ with initial values $\xi(t_{1}+l,t_{1}+l)=1 $ and $\xi (t_{1},t_{1}+r)=0$.\newline
ii) The Hessenberg matrix $\mathbf{\Phi} _{t_{1}+l,t_{1}-r}$, for $r=0,\ldots,t_{1}-t_{2}$ and $r+l\ge 1$, is defined as:
\begin{equation*}
\mathbf{\Phi} _{t_{1}+l,t_{1}-r}=\left[ 
\begin{array}{cc}
\mathbf{\Phi} _{r}^{(2)} & \bar{\mathbf{0}} \\ 
\widetilde{\mathbf{0}} & \mathbf{\Phi} _{l}^{(1)}%
\end{array}\right],
\end{equation*}
where (for $r,l\neq 0$) $\bar{\mathbf{0}}$ is an $r\times l$ matrix of zeros except for $-1$ in its $r\times 1$ entry, and $\widetilde{\mathbf{0}}$ is an $l\times
r$ matrix of zeros except for $\phi _{2,1}$ in its $1\times r$ entry. Notice that $\mathbf{\Phi} _{t_{1}+l,t_{1}-r}$ is a block square matrix of order $l+r$. Its determinant is $\xi(t_{1}+l,t_{1}-r)=\left \vert \mathbf{\Phi}_{t_{1}+l,t_{1}-r}\right \vert $ with initial value $\xi (t_{1},t_{1})=1$.
\end{definition}

{\rm Applying Theorem \ref{TheoGenSol} to the DAB-AR($2;2$) model we obtain the
following Corollary.}

\begin{corollary}
\label{CORGENSOLABAR} The explicit representation of $y_{t_{1}+l}$ in eq. {\rm(\ref{AB-AR(2,2)})} in terms of the prescribed random variables  $y_{t_{2}}$, $y_{t_{2}-1}$, is
given by 
\begin{equation}\label{GENSOLAB_AR(2,2)}
y_{t_{1}+l}=\sum_{r=t_{2}+1}^{t_{1}+l}\xi (t_{1}+l,r)(\varphi(r)+\varepsilon _{r})+\xi (t_{1}+l,t_{2})y_{t_{2}}+\phi _{2,1}\xi(t_{1}+l,t_{2}+1)y_{t_{2}-1}, 
\end{equation}
where $\varphi(r)=\left\{\begin{array}{ll} \varphi_2 & \text{if}\ \  r\le t_1\\
\varphi_1 & \text{if}\ \  r > t_1\end{array}  \right.$.
\end{corollary}

\subsection{Second Moment Structure \label{SecMomStr}}

In this Section we will examine the second moment structure of the DAB-AR ($2;2$) model. To obtain the time-varying variance of $y_{t_{1}+l}$, we will
directly apply Corollary \ref{CORGENSOLABAR}.

First, let $1-\phi _{1,i}B-\phi_{2,i}B^{2}=(1-\lambda _{1,i}B)(1-\lambda_{2,i}B)$, \ for $i=1,2,3$.

\bigskip

\textbf{Assumption 1} (Second-Order Stationarity): $\left \vert \lambda_{m,i}\right\vert <1$, $m=1,2$, for $i=1,3$.\newline
The above condition implies that the DAB-AR($2;2$)-process is second-order stationary.

The following Proposition states expressions for the time-varying variance
of $y_{t_{1}+l}$ in eq. (\ref{GENSOLAB_AR(2,2)}).

\begin{proposition}
Consider the model in eq. {\rm(\ref{AB-AR(2,2)})}. Under Assumption 1, the $\var(y_{t_{1}+l})$ is given by \begin{equation}
\var(y_{t_{1}+l})=A_{t_{1}+l}\sigma _{1}^{2}+B_{t_{1}+l}\sigma
_{2}^{2}+C_{t_{1}+l}\sigma _{3}^{2},  \label{VARDABAR2}
\end{equation}%
where\vspace{-0.15in}
\begin{align*}
A_{t_{1}+l}& =\sum_{r=1}^{l}\xi ^{2}(t_{1}+l,t_{1}+r),\text{ }%
B_{t_{1}+l}=\sum_{r=0}^{t_{1}-t_{2}-1}\xi ^{2}(t_{1}+l,t_{1}-r)\text{,} \\
C_{t_{1}+l}& =\frac{\left[ (1-\phi _{2,3})\left( \xi
^{2}(t_{1}+l,t_{2})+\phi _{2,1}^{2}\xi ^{2}(t_{1}+l,t_{2}+1)\right) +2\phi
_{1,3}\xi (t_{1}+l,t_{2})\phi _{2,1}\xi (t_{1}+l,t_{2}+1)\right] }{(1+\phi
_{2,3})[(1-\phi _{2,3})^{2}-\phi _{1,3}^{2}]}\text{.}
\end{align*}%
Further, if in the above expression we set: $t_{1}=t_{2}$, and therefore $\phi_{m,1}=\phi _{m,2}$ for $m=1,2$, and $\sigma_{1}=\sigma_{2}$, then we obtain the $\var(y_{t_{2}+l})$, which is
equivalent to the case of one break (notice that in this case $B_{t_{2}+l}=0$%
):%
\begin{equation*}
\var(y_{t_{2}+l})=A_{t_{2}+l}\sigma _{2}^{2}+C_{t_{2}+l}\sigma
_{3}^{2}.
\end{equation*}%
Finally, if in addition we set $l=0$ then we obtain the $\var(y_{t_{2}})$, which (since $A_{t_{2}}=0$, $\xi _{t_{2},t_{2}}=1$, $\xi
_{t_{2},t_{2}+1}=0$) is the well known formula for the time invariant AR($2$) model:

\begin{equation*}
\var(y_{t_{2}})=\frac{(1-\phi _{2,3})\sigma _{3}^{2}}{(1+\phi
_{2,3})[(1-\phi _{2,3})^{2}-\phi _{1,3}^{2}}.
\end{equation*}
\end{proposition}

In the next Section we will show how the above results can be used to derive
a time-varying second-order measure of persistence.

\subsection{Time-varying Persistence }\label{TVP copy(1)}

The most often applied time invariant measures of first-order (or mean)
persistence are the largest autoregressive root (LAR), and the sum of the
autoregressive coefficients (SUM); see, e.g., Pivetta and Reis (2007). As
pointed out by Pivetta and Reis in relation to the issue of recidivism by
monetary policy its occurrence depends very much on the model used to test
the natural rate hypothesis, i.e., the hypothesis that the SUM or the LAR
for inflation data is equal to one. Obviously, if both measures ignore the
presence of breaks then they will potentially under or over estimate the
persistence in the levels. The LAR has been used to measure persistence in
the context of testing for the presence of unit roots (see, for details,
Pivetta and Reis, 2007).

In the following, we suggest a time-varying second-order (or variance)
persistence measure that is able to take into account the presence of breaks
not only in the mean but in the variance as well. Fiorentini and Sentana
(1998) argue that any reasonable measure of shock persistence should be
based on the IRF's. For a univariate process $x_{t}$ with \textit{i.i.d.}
errors, $e_{t}$, they define the persistence of a shock $e_{t}$ on $x_{t}$ as $P(x_{t}\left \vert e_{t}\right. )\overset{\mathrm{def}}{=}\displaystyle%
\frac{\var(x_{t})}{\var(e_{t})} $. Clearly $P(x_{t}\left \vert e_{t}\right. )$ will take its minimum value of one if $x_{t}$ is white noise and it will not exist (will be infinite) for an I($1$) process. It follows directly from eq. (\ref{VARDABAR2}) that $P(y_{t_{1}+l}\left \vert \varepsilon_{t_{1}+l}\right.)= \displaystyle \frac{%
\var(y_{t_{1}+l})}{\mathbb{\sigma }_{1}^{2}}$, is given by\vspace{-0.1in}
\begin{equation}\label{P}
P(y_{t_{1}+l}\left \vert \varepsilon _{t_{1}+l}\right.
)=A_{t_{1}+l}+B_{t_{1}+l}\frac{\sigma _{2}^{2}}{\sigma _{1}^{2}}+C_{t_{1}+l}%
\frac{\sigma _{3}^{2}}{\sigma _{1}^{2}}. \vspace{-0.1in} 
\end{equation}
If Assumption 1 is violated then conditional measures of second-order persistence can be constructed using the variance of the forecast error instead of the unconditional variance (results not reported but are available upon request).

Having derived explicit formulas for time-varying second-order (or variance)
persistence measures, in the next Section we show the empirical relevance of these results using U.S. inflation data.\footnote{Cogley and Sargent (2002) measured persistence by the spectrum at frequency zero, $S_{0}$. As an example, for the time invariant AR($2$) model this will
be given by: $S_{0}=\displaystyle \frac{\sigma _{\varepsilon }^{2}}{2\pi(1-\phi _{1}-\phi _{2})^{2}}$.}

\section{ Inflation Data}

\label{sec. InflationData}

In this Section we directly link econometric theory with empirical evidence.
In our empirical application we consider the possible presence of structural
breaks in inflation for the United States. We use quarterly data on the GDP
deflator as the measure of price level. The data set consists of
observations from 1963Q4 to 2018Q1. Inflation is calculated as the quarterly
change of price level at an annualized rate calculated as $\pi _{t}=400(\ln
(P_{t}/P_{t-1}).$

In terms of inflation modelling, the period under consideration is of
particular interest as it covers the boom-time inflation of the late 1960s,
the stagflation in the 1970s, and the double-digit inflation of the early
1980s. During this period substantial shifts in monetary policy occurred,
most notably the Fed's radical step of switching policy from targeting
interest rates to targeting the money supply in the early 1980s. Therefore
when modelling inflation it is important to allow for time-varying
parameters.

\subsection{Unit Root Tests}
Although we allow for regime shifts, we are particularly interested in
modelling changes in inflation persistence. In the related literature
inflation persistence is defined as the tendency of inflation to converge at
a slow pace to the long run equilibrium level after a shock. Monetary policy
authorities are particularly interested in knowing the speed at which the
inflation rate converges to the central bank's inflation target following
macroeconomics shocks. However, as shown in Levin and Piger (2004) not
accounting for structural breaks may lead to overestimating inflation
persistence. In this regards, a growing body of research has found evidence
that the monetary policy target has an impact on the persistence properties
of inflation as well as on its volatility (see for example Brainard and
Perry, 2000 or Taylor, 2000). That monetary policy actions affect
persistence of inflation is of interest as it has important implications for
inflation modelling as changes in regimes of monetary policy may leave
econometric models open to the Lucas critique.

In the empirical literature a common approach for modelling inflation persistence is to estimate a univariate AR$(p)$ model where the sum of the
estimated autoregressive coefficients is used to approximate the sluggishness with which the inflation process responds to macroeconomic shocks (see for
example Pivetta and Reis, 2007) and/or apply unit root tests.

In Table \ref{table: Unit root} a number of common unit roots tests are reported. Namely: ADF (Augmented Dickey--Fuller), ERS, by Elliott \textit{et al}. (1996), and MZ GLS, suggested by Perron and Ng (1996) and Ng and Perron (2001). As recommended by Ng and Perron (2001), the choice of the number of lags is based on the modified Akaike information criterion (AIC). 
 \renewcommand{\arraystretch}{1.2} 
\begin{table}[h!]
	\centering
	\begin{tabular}{lclccc}
		\multicolumn{6}{l}{\textbf{Table \ref{table: Unit root}}.\ Unit root Tests.} \vspace{0.01in}\\ 
		\hline\hline
		Test Statistic & \multicolumn{5}{c}{} \\
		\cline{1-1}
		\mbox{}\vspace{-0.15in}\\
		ADF &  & $-3.229^{\ast \ast }$ &  &  &  \\ 
		ERS &  & $-3.154^{\ast \ast }$ &  &  &  \\ 
		MZ$_{a}$ &  & $-19.642^{\ast }$ &  &  &  \\ 
		MZ$_{t}$ &  & $-3.1331^{\ast }$ &  &  &  \\ 
\hline\hline
\end{tabular}
\caption{\footnotesize The notations $^*, ^{**}$ indicate the statistical significance at 1\% and 5\%, respectively.}
\label{table: Unit root}\vspace{-0.1in}
\end{table}

Table \ref{table: Unit root} shows that, in general, we can reject the null hypothesis of a unit root in inflation series.

\subsection{ Structural Breaks}

In what follows we estimate AR($p$) models with abrupt structural breaks.
The optimal model is the DAB-AR($2;2$) model in eq. (\ref{AB-AR(2,2)}). The
choice of the number of lags was based on the modified AIC and the Bayesian
information criteria. The break points are treated as unknown. Note that
breaks in the variance are permitted provided that they occur at the same
dates as the break in the autoregressive parameters. Benati (2008) also used
an AR model allowing for time-varying volatility. Cogley and Sargent (2005)
also estimated a model in which the variance of innovations can vary over
time. For each $l$ partition $\left( T_{1},...,T_{l}\right) $ the DAB-AR($2;l$) model can be estimated using the least-squared principle by minimizing the sum of the squared residuals where the minimization is taken over all
partitions. Since the break points are discrete parameters and can only take
a finite number of values they can be estimated by grid search using dynamic
programming (see Bai and Perron, 2003, for more details).

Coming to the estimation procedure, the first step is to identify possible
points of parameter changes. In order to do so the Bai and Perron (2003)
sequential tests on inflation rates is used to identify possible breaks
during the sample period. They propose an $F$-type\ test for $l$ versus $l+1$
breaks, which we refer to as $\sup F_{t}(l+1|l)$. The testing procedure
allows for a specific to general modelling strategy for the determination of
the number of breaks in each series. The test is applied to each segment
containing the $T_{i-1}$ to $T_{i}$ $(i=1,...,l+1)$. In particular, the
procedure involves using a sequence of $(l+1)$ tests, where the conclusion
of a rejection in favour of a model with $(l+1)$ breaks if the overall
minimal value of the sum of squared residuals is sufficiently smaller than
the sum of the squared residuals from the $l$ break model.
Note that the sum of the squared residuals is calculated over all segments where an additional break is included and compared with the residuals from the $l$ model. Therefore, the break date selected is the one associated with the overall minimum.

The results of the structural break test are summarized in Panel A of Table \ref{table: structural break} below. The first column in Panel A compares the null hypothesis of $l$ breaks against the alternative hypothesis of $l+1$ breaks, the second column reports the calculated value of the statistics and the third column the critical value of the test. 
\begin{table}[h!]
	\centering		
\begin{tabular}{cccccc}
\multicolumn{6}{l}{\textbf{Table \ref{table: structural break}}.\ Structural break test and estimation results.} \\ 
		\hline\hline
\multicolumn{6}{l}{\textbf{Panel A:} Bai and Perron tests of $L+1$\ vs. $L$\ sequentially determined breaks} \\ 
Null hypotheses &  & F-Statistic & Critical Value &  &  \\ 
\cline{1-1}\cline{3-6}
&  &  &  &  &  \vspace{-0.15in}\\ 
$H_{0}:0$\ $\rm vs$\ $1$ &  & $57.96^{\ast \ast }$ & 13.98 &  &  \\ 
$H_{0}:1$\ $\rm vs$\ $2$ &  & $18.13^{\ast \ast }$ & 15.72 &  &  \\ 
$H_{0}:2$\ $\rm vs$\ $3$ &  &\hspace{-0.17in} $13.57$ & 16.83 &  & \\
\hline 
\multicolumn{6}{l}{\textbf{Panel B:} Model Estimation and Misspecification Tests} \\ 
Period &  & $\varphi _{i}$ & $\phi _{1,i}$ & $\phi _{2,i}$ & $\sigma _{i}$\\ 
\cline{1-1}\cline{3-6}
\multicolumn{1}{c}{} &  &  &  &  &  \vspace{-0.15in}\\ 
\multicolumn{1}{c}{1964Q2-1976Q3} &  & $\underset{(0.224)}{0.496^{\ast\ast }}$ & $\underset{(0.108)}{0.470}^{\ast }$ & $\underset{(0.102)}{0.376^{\ast }}$ & $\underset{(0.367)}{1.077^{\ast} }$ \\ 
\multicolumn{1}{c}{1976Q4-1986Q2} &  & $\underset{(0.954)}{3.637^{\ast }}$ & 
$\underset{(0.119)}{0.710^{\ast }}$ & $\underset{(0.112)}{0.127}$ & $\underset{(0.689)}{2.300^{\ast }}$ \\ 
\multicolumn{1}{c}{1986Q3-2018Q1} &  & $\underset{(0.396)}{2.859^{\ast }}$ & 
$\underset{(0.082)}{0.247^{\ast }}$ & $\underset{(0.077)}{-0.314^{\ast }}$ & $\underset{(0.489)}{2.160^{\ast }}$ \\ 
\hline
\multicolumn{2}{c}{$R^{2}$} & \multicolumn{1}{c}{} & $0.614$ &  &  \\ 
\multicolumn{2}{c}{Breusch-Godfrey Test} & \multicolumn{1}{c}{} & $\underset{(0.561)}{2.055}$ &  &  \\ 
\multicolumn{2}{c}{White Test} & \multicolumn{1}{c}{} & $\underset{(0.376)}{3.103^{{}}}$ &  &  \\ \hline\hline
\end{tabular}
\caption{\footnotesize \textbf{Panel A} reports the calculated Bai-Perron test for structural
breaks along with the critical value of the test taken from Bai and Perron (2003).  \textbf{Panel B} provides the estimated parameters along with the associated standard errors (see the next Subsection). The notations $^*, ^{**}$ indicate the statistical  significance at 1\% and 5\%,respectively. The $p$ values for the misspecification tests are given in parenthesis.}
\label{table: structural break}\vspace{-0.1in}
\end{table}
\newpage
Observing the calculated values of the test, it appears that the null hypothesis zero versus one break is rejected in favour of the alternative hypothesis. Similarly, the hypothesis of one break versus two breaks is rejected. However, the null hypothesis of two versus three breaks in not rejected, therefore we conclude that there are two structural breaks. The first break occurred in the mid-1970's, when the Fed tightened monetary policy to fight the high inflation rate after the end of the Bretton Woods period. The second break occurred in 1986 when the Fed embarked on an aggressive policy to reduce inflation, which reached unusually high levels starting from the 70s.

As a conclusion, inflation fell from $10.5\%$ at the end of 1980 to $1.1\%$ in 1986Q2, which is also the date of the estimated break.\footnote{McConnell and Perez-Quiros (2000) have detected a fall in the volatility of output \ after 1984 as well.}

\subsection{ Estimation Results}
Panel B of Table \ref{table: structural break} reports  
the estimated model parameters along with the relative misspecification tests. According to the reported  parameter estimates, the inflation process is well approximated by a second-order autoregression. Moreover, the drift parameters $\varphi_{i}$, $i=1,2,3$ increase from $\varphi_{3}=0.496$ before 1976Q3 to $\varphi_{2}=3.637$ during the period 1976-1986. The increase in the drift reflects the fact that toward the second half of the 70s until the middle of the 80s the inflation level was stubbornly high. After 1986 the smaller magnitude of the estimated drift reflects the lower average inflation rates that the US enjoyed over the last three decades. This is in line with the finding in Levin and Piger (2004), who provide statistical evidence for a fall in the intercept after the early 1990s. Kozicki and Tinsley (2002) interpreted this shift as change in the long-run inflation target of the Federal Reserve.

Considering now the estimated autoregressive parameters, $\phi _{1,i}$ and $\phi _{2,i}$, according to the estimates until 1986 the inflation process
had a high \textit{intrinsic} persistence ($\phi _{1,3}+\phi_{2,3}=0.846\simeq \phi _{1,2}+\phi _{2,2}=0.837$), but it has fallen ever
since. These results are consistent with the findings in Cogley and Sargent (2002) (see also Brainard and Perry, 2000, and Taylor, 2000). With respect
to the variance parameter $\sigma _{i}$, we see that the volatility of the innovation was relatively high during the decade 1976-1986 ($\sigma
_{2}=2.30 $) and it has fallen slightly in the last thirty years ($\sigma _{3}=2.160$). However, it did not go back to the relatively low level before
1976 ($\sigma _{1}=1.077$). This is probably due to the fact that the last period included the turmoil of the financial crisis that started in 2005
(see, for example, Stock and Watson, 2009).

Our estimated model confirms that changes in inflation dynamics can be explained by changes in the drift, the \textit{intrinsic} persistence and
the variance parameter. To summarize our results, we find evidence that the parameters in the models capturing persistence change over time. Therefore,
not allowing for time-varying coefficients in the estimation procedure would result in a less accurate modelling of the inflation process. This, in light of the simulation results in Section \ref{Forecasting}, may lead to poor forecasting.
Finally, the misspecification tests are reported at the bottom of Panel B.
It appears that the Breusch-Godfrey for autocorrelation does not reject the
null hypothesis of no serial correlation. Similarly, the White test for
heteroscedasticity does not reject the null hypothesis of homoscedasticity,
therefore indicating that the model does not suffer from misspecification.

\subsection{Forecasting}
\label{Forecasting}
We now consider the out-of-sample forecasting performance of the estimated model
(see Table \ref{table: structural break}). In order to investigate the effect of model
misspecification on the forecasted inflation level we compare three models.
The first model, which we label as Model 1, is the estimated DAB-AR(2;2).
The second model, which we refer to as Model 2, is the true model which we obtained by simulating the inflation process using the estimated parameters
in Table \ref{table: structural break} as a data generating process. Finally, the third model is the
misspecified AR(2) model with no time-varying parameters, which we label as Model 3.

The evaluation of the out-of-sample forecast exercise does not rely on a single criterion; for robustness we compare the results of three different forecasting measures, namely, the root mean square forecast error (RMSE), the mean absolute error (MAE) and Theil's Inequality Coefficient (U Coefficient). 

Table \ref{table: forecasting inflation} below reports the results of the forecasting exercise.\footnote{For forecasting under structural breaks, see for example, Pesaran and Timmermann (2005).} In columns 1 and 2 the forecasting horizon and the performance measure are shown, respectively. Columns 3-5 show the forecasting results.
It follows directly form the results of Table \ref{table: forecasting inflation} that according to the RMSE and MAE criteria the DAB-AR (2;2) model performs better than its misspecified counterpart.
According to these two performance measures Model 1 has forecasting properties in line with those obtained using the true model, Model 2.
However, looking at the U coefficient measure the results are more mixed with Model 3 outperforming Model 1 in the short horizon and Model 1 having
superior performance in the longer horizon.

\begin{table}[h!]
\centering	
\begin{tabular}{ccccccc}
\multicolumn{7}{l}{\textbf{Table \ref{table: forecasting inflation}.}  Forecasting inflation in the United States: point predictive performances.} \\ \hline\hline
{\small Forecast Horizon} & {\small Forecast Error Measure} & {\small Model 1} &  & {\small Model 2} &  & {\small Model 3} \\ 
\cline{1-1}\cline{2-3}\cline{5-5}\cline{7-7}
&  &  &  &  &  &  \vspace{-0.15in}\\ 
{\small 1} & {\small RMSE} & {\small 0.0134} &  & {\small 0.0110} &  & 
{\small 0.0194} \\ 
{\small 4} &  & {\small 0.0141} &  & {\small 0.0121} &  & {\small 0.0167} \\ 
{\small 8} &  & {\small 0.0132} &  & {\small 0.0149} &  & {\small 0.0242} \\ 
&  &  &  &  &  &  \vspace{-0.1in}\\ 
{\small 1} & {\small MAE} & {\small 0.0166} &  & {\small 0.0111} &  & 
{\small 0.0944} \\ 
{\small 4} &  & {\small 0.0112} &  & {\small 0.0101} &  & {\small 0.0144} \\ 
{\small 8} &  & {\small 0.0112} &  & {\small 0.0104} &  & {\small 0.0208} \\ 
&  &  &  &  &  &  \vspace{-0.1in}\\ 
{\small 1} & {\small U Coefficient } & {\small 0.323} &  & {\small 0.251} & 
& {\small 0.293} \\ 
{\small 4} &  & {\small 0.258} &  & {\small 0.241} &  & {\small 0.243} \\ 
{\small 8} &  & {\small 0.327} &  & {\small 0.287} &  & {\small 0.407} \\ 
\hline\hline
\end{tabular}
\caption{\footnotesize The table compares the out-of-sample point forecasts of three models.\textbf{ Model 1} is the DAB-AR (2;2)  estimated model (see Table 2).\textbf{ Model 2} is obtained using simulated data. \textbf{ Model 3} is an AR(2) process with no time-varying parameters. The forecast measures are $i)$  the RMSE, $ii)$ the MAE, and  $iii)$\ the U Coefficient. The forecast horizon is $1,4,$ 	and $8$\ quarters.}
\label{table: forecasting inflation}
\end{table}
\newpage
Having investigated the out-of-sample forecasting performance of the DAB-AR-(2;2) model we next investigate whether inflation and its volatility
are highly persistent.

\subsection{Inflation Persistence}

Pivetta and Reis (2007) employ different estimation methods and measures of
persistence. Estimating the persistence of inflation over time using
different measures and procedures is beyond the scope of this paper.\footnote{Pivetta and Reis (2007) applied a Bayesian approach, which explicitly treats
the autoregressive parameters as being stochastically varying and it
provides their posterior densities at all points in time. From these, they
obtained posterior densities for the measures of inflation persistence. Such
estimates of persistence are forward-looking, since they are meant to
capture the perspective of a policy maker who at a point in time is trying
to foresee what the persistence of inflation will be. They also estimated
-looking measures of persistence that the applied economist forms at a point in time, given all the sample until then.
Pivetta and Reis (2007) also used an alternative set of estimation
techniques for persistence. They assumed time invariant autoregressive
parameters and re-estimated their AR model on different sub-samples of the
data, obtaining median unbiased estimates of persistence for each
regression. Finally, Pivetta and Reis also employed rolling and recursive
unit root tests.} In this Section we depart from their study in an important
way, that is we contribute to the measurement over time of inflation
persistence by taking a different approach to the problem and estimate a
DAB-AR model of inflation dynamics grounded on econometric theory, and we
compute an alternative measure of persistence, that is, the second-order
persistence (using the methodology in Sections \ref{SecMomStr} and \ref{TVP
copy(1)}), which distinguishes between changes in the dynamics of inflation
and its volatility (and their persistence).

As pointed out in the above cited reference estimates of the inflation persistence affect the tests of the natural hypothesis neutrality. Therefore
detecting whether persistence has recently fallen is key in assessing the
likelihood of recidivism by the central bank. In addition, if the central
bank feels encouraged to exploit an illusory inflation-output trade off, the
result could be high inflation without any accompanying output gains.
Furthermore, research on dynamic price adjustment has emphasized the need
for theories that generate inflation persistence.

Table 4 presents the within each period time invariant first and
second-order measures of persistence for all three periods. The first three
columns report the three first-order measures of persistence (LAR, $%
1/(1-SUM) $ and $\expectation(\pi_{t})$). For the first two measures Period 1
yields the highest persistence. In particular, the persistence (measured by $%
1/(1-SUM)$) decreases by $5.5\%$ in the post-1976 period and it decreases
further by $85\%$ in the post-1986 period. The mean of inflation, $\expectation%
(\pi_{t})$, increases by $59.3\%$ in the second period and it decreases by $%
88\%$ in the third period. The last three columns of Table 4 report the
three second-order measures of persistence, i.e. $S_{0}$, $%
P_2(\pi_{t}\left
\vert \varepsilon_{t}\right.)\overset{\mathrm{def}}{=}%
\displaystyle \frac{\var(\pi_{t})}{\var(\varepsilon_{t})}$
and $\var(\pi_{t})$. For two out of the three measures the post-1986 period exhibits the lowest persistence whereas in the second period the
persistence is the highest. The variance of inflation, $Var(\pi _{t})$, from 1976 to 1986 is almost five times the variance of inflation of the pre-1976 period and it is almost three times the variance of the post-1986 period.

\begin{table}[h!]
	\centering	 
\begin{tabular}{c|ccc|ccc}
	\multicolumn{7}{l}{\textbf{Table \ref{persistence}.}  Persistence for each of the three 		periods/models.} \\ 
	\hline\hline
\multicolumn{7}{c}{First and Second-order Measures of Persistence} \\ 
\hline
\multicolumn{1}{c|}{Period}	& \multicolumn{3}{c|}{First-Order} & \multicolumn{3}{c}{Second-Order} \\
\hline 
\multicolumn{1}{c|}{} & LAR & 1/(1-SUM) & $\expectation(\pi_{t})$ & $S_{0}$ & $P(\pi_{t}\left \vert \varepsilon _{t}\right. )$ & $\var(\pi_{t})$ \\ 
$1964Q_{2}-1976Q_{3}$ & 0.892 & 6.493 &\hspace{0.05in} 3.221 & \hspace{0.07in}7.784 & 2.692 & \hspace{0.03in}3.122 \\ 
$1976Q_{4}-1986Q_{2}$ & 0.858 & 6.135 & 22.313 & 31.688 & 3.002 & 15.881 \\ 
	$1986Q_{3}-2018Q_{1}$ & 0.560 & 0.937 &\hspace{0.05in} 2.679 & \hspace{0.08in}0.652 & \hspace{0.01in}1.150 & \hspace{0.02in} 5.365 \\ 
	\hline\hline
\end{tabular}
\caption{\footnotesize For each period, $n=1,2,3$  we use the six alternative measures to calculate the (within each period time invariant) first and second-order persistence.}
\label{persistence}
\end{table}

The following graphs\footnote{The graphs have been designed and plotted using the Mathematica drawing tools.} 
of the measures $P_2(\pi_{t}\left \vert \varepsilon_{t}\right.)$ and $\var(\pi _{t})$ reflect the dynamics of the second-order time-varying inflation persistence. In the x-axis each unit represents a year-quarter starting with $1964Q_{2}$, chosen as the first. In particular, $1976Q_{3}=49$ (49-th year-quarter) and $1986Q_{2}=88$ (88-th year-quarter).
\begin{figure}[h]
\begin{subfigure}[h]{0.45\linewidth}
		\includegraphics[width=\linewidth]{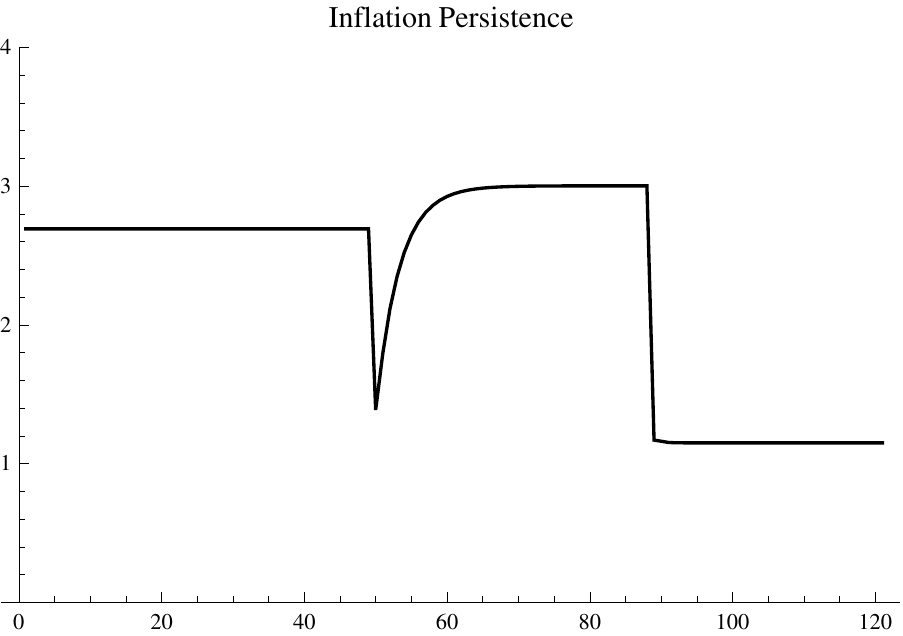}
		\caption{$P_2(\pi_{t}\left \vert \protect \varepsilon_{t}\right. )$}
	\end{subfigure}
\hfill 
\begin{subfigure}[h]{0.45\linewidth}
		\includegraphics[width=\linewidth]{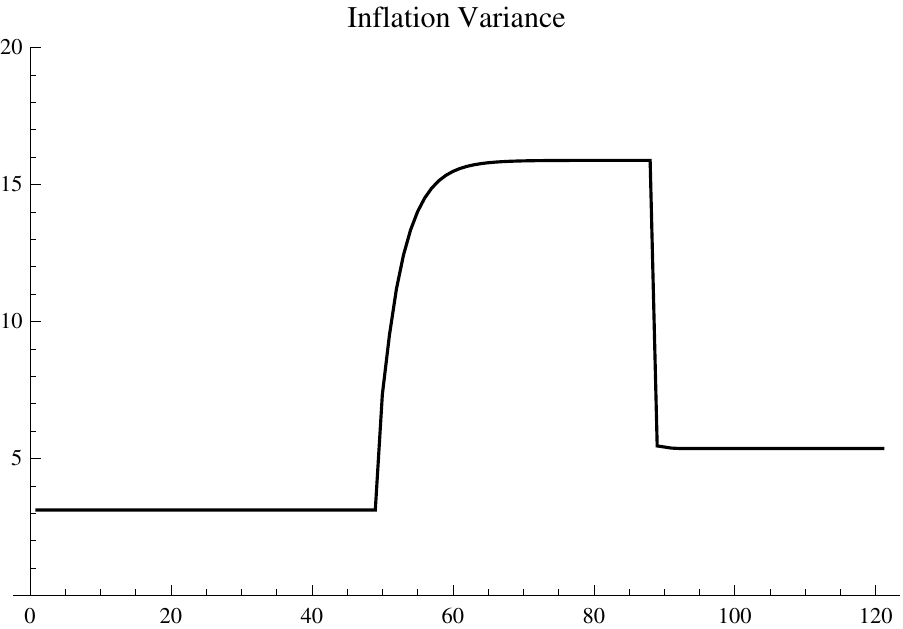}  
		\caption{$Var(\pi_{t})$}
	\end{subfigure}
\caption{Second-Order time-varying Persistence}
\end{figure}
Details of how we construct the graphs are presented in the Online Appendix K.

\bigskip 

Some key features for the graph of the inflation variance $Var(\pi_{t})$ are discussed below: i) In the pre-76 period the graph is constant: $Var(\pi_{t})=3.122$. ii) Within the post-76 and pre-86 period, the graph
increases abruptly next to the quarter $1976Q_4$, but at a decreasing rate, reaching in the end the highest value $Var(\pi_{t})=15.881$. iii) In the post-86 period the graph stabilizes to $Var(\pi_{t})=5.365$, after an abrupt
drop next to the quarter $1986Q_3$. Analogous statements can be addressed for the inflation persistence graph $P_2(\pi_{t}\left \vert \varepsilon_{t}\right.)$. As illustrated above, the main difference between
the shapes of the two graphs is due to the abrupt drop next to the quarter $1976Q_4$ followed shortly afterwards by an abrupt increase at a decreasing rate.

The graphs of the two measures $P_1(\pi_{t}\left \vert\varepsilon_{t}\right.)%
\overset{\mathrm{def}}{=}\displaystyle \frac{\expectation(\pi_{t})}{\varphi(t)}
$ and $\expectation(\pi_{t})$ for the first-order persistence are shown below:
\begin{figure}[h]
\begin{subfigure}[h]{0.45\linewidth}
		\includegraphics[width=\linewidth]{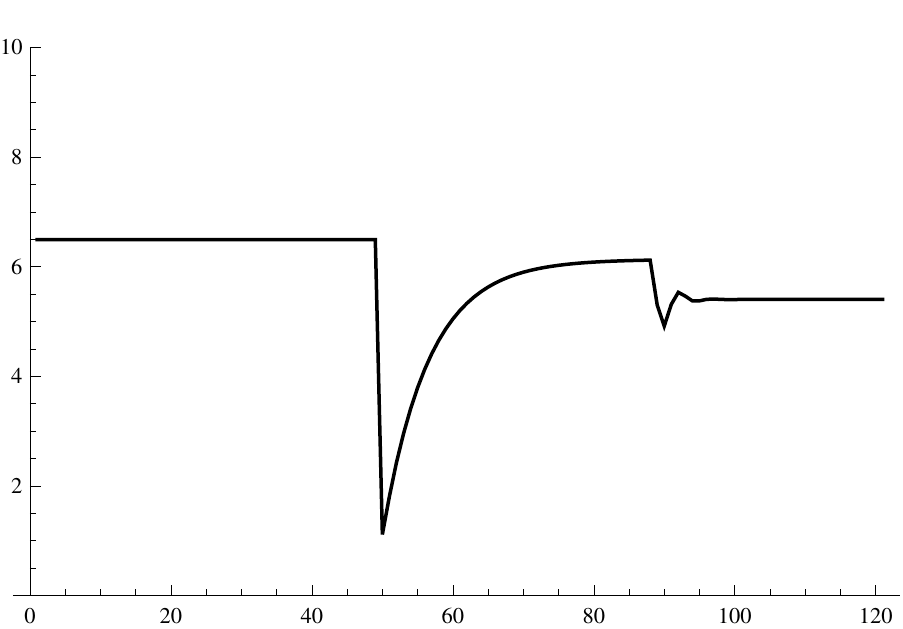}
		\caption{$P_1(\pi_{t}\left \vert \protect \varepsilon_{t}\right. )$}
	\end{subfigure}
\hfill 
\begin{subfigure}[h]{0.45\linewidth}
		\includegraphics[width=\linewidth]{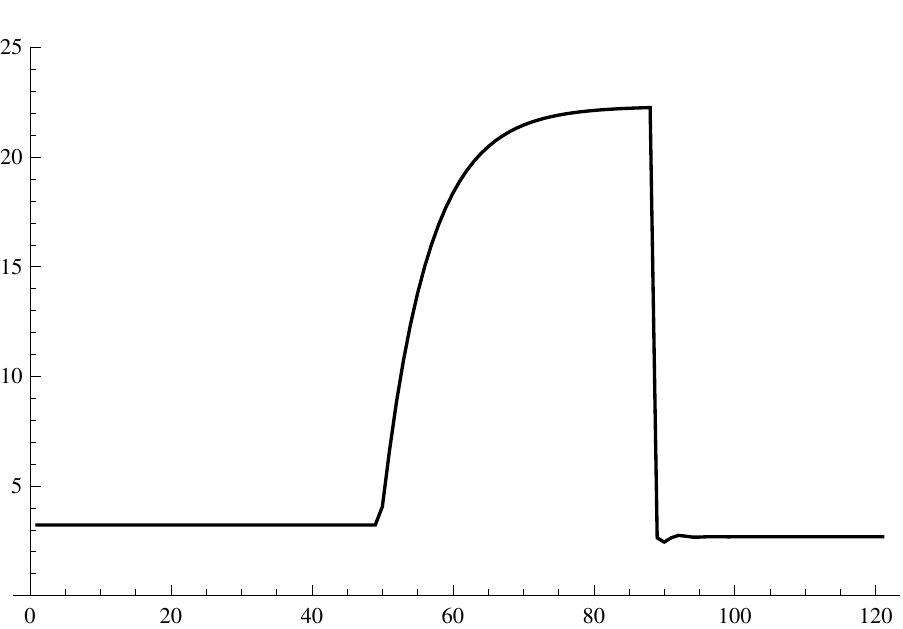}  
		\caption{$\expectation(\pi_{t})$}
	\end{subfigure}
\caption{First-Order time-varying Persistence}
\end{figure}

In sum, our main conclusion is that for our chosen specification (DAB-AR
model) the preferred measure of persistence, that is the second-order
persistence, as measured by the conditional variance of inflation, increased
considerably from 1976 onwards, whereas in the post-1986 period the
persistence falls to even lower levels than the pre-1976 period. Our results
are in line with those in Cogley and Sargent (2002), who find that the
persistence of inflation in the United States rose in the 1970s and remained
high during this decade, before starting a gradual decline from the 1980s
until the early 2000s (similar to the results of Brainard and Perry, 2000,
and Taylor, 2000). Stock and Watson (2002) found no evidence of a change in
persistence in U.S. inflation. However, they found strong evidence of a fall
in volatility. Therefore their results are in agreement with ours.

\section{ Conclusions and Future Work\label{SECCONCL}}
It is important to understand the fundamental properties of  time series linear models with variable coefficients in order to treat them  effectively, and other more complicated structures originated from them. We have put forward a methodology for solving linear stochastic time-varying difference equations. The theory presented makes no claim to being applicable in all linear processes with variable coefficients. However, the cases covered are those which belong to the large family of ``time-varying" models with ARMA representations. Our methodology is a practical tool that can be applied to many dynamic problems. As an illustration we studied an AR specification with abrupt breaks, which is grounded on econometric theory. The second moment structure of this construction was employed to obtain a new time-varying measure of second-order persistence.

To summarize, we identified a lack of a universally applicable approach yielding an explicit solution to TV-ARMA models. Our response was to try and fill the gap by developing a coherent body of theory, which implicitly contains the invertibility of a time-varying polynomial, and, therefore, can replace the convenient tool of characteristic polynomials. In particular, the unified theory does the following things: First, it provides a new technique that gives an explicit solution representation  of such models for any set of prescribed random variables $\{y_r\}_{s+1-p\le r\le s}$. This enables us to treat within a unified scheme of ARMA processes with variable coefficients including the cases of deterministic and stochastic ones. Second, it derives a sufficient condition for their asymptotic stability  and highlights its strong connection with the Wold-Cram\'{e}r solution decomposition. Third it generates explicitly the second moments of these processes along with necessary and sufficient conditions for their existence, which (in the case of the deterministically varying coefficients) are required for the quasi maximum likelihood and central least squares estimation. Fourth it provides conditions for the invertibility of these processes followed by optimal forecasts based on a finite set of observations. Fifth, it establishes a sufficient condition for their asymptotic efficiency grounded on the boundedness of the mean square error.

We developed this new technique, which can be applied virtually unchanged in
every ``ARMA" environment, that is to the even larger family of ``time
varying" models, with ARMA representations (i.e., GARCH type of [or
stochastic] volatility and Markov switching processes; for the abundant
literature on weak ARMA representations see, for example, Francq and Zako%
\"{\i}an, 2005, and the references therein). Thus our results are applied to
TV-GARCH models as well without any significant difficulties. This generic
framework that forms a base for such a general approach releases us from the
need to work with characteristic polynomials and, by enabling us to examine
a variety of specifications and solve a number of problems, helps us to
deepen our familiarity with their distinctive features.

The empirical relevance of the theory has been illustrated through an
application to inflation rates. Our estimation results led to the conclusion
that U.S. inflation persistence has been high since 1976, whereas after 1986
the persistence falls to even lower levels than the pre-1976 period, a
finding which agrees with those of Brainard and Perry (2000), Taylor (2000)
and Cogley and Sargent (2002).

The usefulness of our unified theory is apparent from the fact that it
enables us to analyze an abundance of models and solve a plethora of
problems. In addition, an extension of the methodology developed in this
paper enables us to (just to mention a few consequences): i) examine in
depth infinite order autoregressions with either constant or variable
coefficients, since it releases us from the need to work with characteristic
polynomials and ii) obtain the fourth moments of TV-GARCH models, which
themselves follow linear time-varying difference equations of infinite
order, taking advantage of the fact that the various GARCH formulations have
weak ARMA representations (see, for example, Karanasos, 1999, and Francq and
Zako\"{\i}an, 2005), iii) work out the fundamental time series properties of
time-varying linear VAR systems (since it can be easily modified and
applied to a multivariate setting; see, for example, Karanasos et al.,
2014), iv) derive explicit formulas for the nonnegativity constraints and
the second moment structure of both constant and time-varying multivariate
GARCH processes (thus extending the results in He and Ter\"{a}svirta, 2004,
Conrad and Karanasos, 2010, and Karanasos and Hu, 2017).

Hallin (1986) applied recurrences in a multivariate context to obtain the
Green's matrices. Work is at present continuing on the multivariate case.
When this has been completed one should be able to apply, without any major
alterations, the methods of this paper to multivariate TV ARMA and GARCH
models. Spectral factorization is another important problem that can be
solved by our new representations.

Some of these research issues are already work in progress and the rest will
be addressed in future work.

\end{document}